\documentclass[11pt,a4paper]{article}
\usepackage[utf8]{inputenc}
\usepackage[T1]{fontenc}
\usepackage{amsmath,amssymb,amsthm}
\usepackage{graphicx}
\usepackage{booktabs}
\usepackage{geometry}
\usepackage{hyperref}
\usepackage{natbib}
\usepackage{subcaption}
\usepackage{multirow}
\usepackage{tabularx}
\usepackage{float}
\usepackage{appendix}
\usepackage{lmodern}

\newcommand{\D}{\mathrm{D}}
\newcommand{\N}{\mathrm{ND}}

\title{A Quantitative Framework to Predict Wait-Time Impacts Due to AI-Triage Devices in a Multi-AI, Multi-Disease Workflow}

\author{Michelle Mastrianni$^{1}$, Rucha Deshpande$^{1}$,Frank W. Samuelson$^{1}$,\\ Yee Lam Elim Thompson$^{1*}$\\ \\
$^{1}$ The U.S. Food and Drug Administration\\
$^{*}$ Correspondence: YeeLamElim.Thompson@fda.hhs.gov}

\date{\today}

\begin{document}

\maketitle

\begin{abstract}
The deployment of multiple AI-triage devices in radiology departments has grown rapidly, yet the cumulative impact on patient wait-times across different disease conditions remains poorly understood. This research develops a comprehensive mathematical and simulation framework to quantify wait-time trade-offs when multiple AI-triage devices operate simultaneously in a clinical workflow. We created multi-QuCAD, a software tool that models complex multi-AI, multi-disease scenarios using queueing theory principles, incorporating realistic clinical parameters including disease prevalence rates, radiologist reading times, and AI performance characteristics from FDA-cleared devices. The framework was verified through four experimental scenarios ranging from simple two-disease workflows to complex nine-disease systems, comparing preemptive versus non-preemptive scheduling disciplines and priority versus hierarchical triage protocols. Analysis of brain imaging workflows demonstrated that while AI-triage devices significantly reduce wait-times for target conditions, they can substantially delay diagnosis of non-targeted, yet urgent conditions. The study revealed that hierarchical protocol generally provides more wait-time savings for the highest-priority conditions compared to the priority protocol, though at the expense of more delays to lower-priority patients with other time-sensitive conditions. The quantitative framework presented provides essential insights for orchestrating multi-AI deployments to maximize overall patient time-saving benefits while minimizing unintended delay for other important patient populations.
\end{abstract}


\section{Introduction}

Artificial intelligence (AI) algorithms are transforming medical imaging through emerging applications such as automated prioritization of patient images. Over the past decade, the U.S. Food and Drug Administration (FDA) has cleared more than 87 AI-enabled computer-aided triage and notification (AI-triage) devices. These devices are intended to identify and prioritize radiological medical images flagged as suspicious by the AI, enabling expedited review of radiological images with time-sensitive disease conditions and potentially improving patient outcome.

Most AI-triage devices are designed to identify a single disease condition and prioritize patient images suspicious of that target disease. To evaluate the wait-time-saving benefits for patients diagnosed with the disease, a mathematical framework based on queueing theory was developed~\cite{Thompson}. This framework considers a single AI-triage device targeting one disease condition and encompasses various clinical workflow factors, such as the traffic intensity, the number of radiologists available for image interpretation, interruptions caused by emergency patients, and differences in radiologist reading times between diseased and non-diseased images. Validated using real-world clinical data~\cite{ThompsonJACR}, the QuCAD software \cite{ThompsonQuCAD} was developed to simulate patient flow dynamics and provide estimates of wait-time savings under these scenarios, assuming a preemptive-resume scheduling discipline in which radiologists are interrupted to address higher-priority cases upon their arrivals.

Previous work primarily focuses on the time-saving benefits for patients diagnosed with the target disease condition when a single AI-triage device is deployed and neglects the potential time delays experienced by patients with other urgent conditions. This time delay becomes even more concerning as clinics have started to integrate multiple AI-triage devices into their workflows. The deployment of multiple devices introduces complex dynamics, necessitating a deeper understanding of the trade-offs between the benefits of reduced wait-times for some patients and the added delays for others. Addressing this balance is essential for optimizing overall patient outcomes and ensuring equitable access to timely care. Moreover, the existing QuCAD software assumes a preemptive-resume discipline. In practice, however, radiologists may adopt a non-preemptive discipline, where they complete the review of the lower-priority case in hand before attending to the newly arrived higher-priority case. These complexities underscore the need for a more realistic model to provide insights on the impact of AI-triage deployment on patient wait-time in a complex clinical workflow. 

To predict the clinical impact of AI-triage devices in complex workflows, this paper presents an extended framework, multi-QuCAD \cite{ThompsonMultiQuCAD}, that provides both simulated and theoretical estimates of wait-time savings and delays across different patient subgroups. Building on the workflow factors in the original QuCAD software, multi-QuCAD expands its scope to account for patients with different disease conditions and the deployment of multiple AI-triage devices in a workflow. Multi-QuCAD also models workflows under both preemptive-resume and non-preemptive scheduling disciplines. To further enhance its applicability, multi-QuCAD simulates two protocol settings for each scheduling discipline. In the priority setting, all AI-positive cases are prioritized to the front of the queue. In contrast, in the hierarchical setting, disease conditions are ranked by their time-sensitiveness, and cases flagged by devices targeting more time-sensitive conditions are reviewed by the radiologists first. 

The main aim of this work is to introduce multi-QuCAD, a quantitative framework to predict the potential wait-time impacts for patients with different conditions when integrating multiple AI-triage devices into radiologist workflows. Our first use case involves stroke, a broad term for a condition that disrupts blood flow to the brain.  Depending on cause, strokes are subclassified into hemorrhagic stroke (either subdural hemorrhage (SDH) or subarachnoid hemorrhage (SAH)) or  large vessel occlusion (LVO).  Each scenario may have one or more AI-triage device(s) targeting each stroke condition. For each scenario, we reported the average time-savings and delays experienced by patients with different stroke types, as well as the agreement between simulation and theoretical predictions. To further assess the consistency between simulation and theory in more complex workflows, we extended our analysis to a hypothetical scenario with nine disease conditions and four AI-triage devices. Multi-QuCAD provides insights into how workflow-specific factors influence patient wait-time, serving as a practical tool to help clinicians optimize the integration of multiple AI-triage devices into their unique clinical workflows, ultimately improving patient care and operational efficiency. 

\section{Problem Statement}
\label{sec:parameters}

Our goal is to predict the wait-time-savings and delay of patient images when one or more AI-triage devices for various disease conditions are deployed within the same workflow. The impact to wait-time-savings and delay will also be explored among different workflow types defined in this section. 

In this work, wait-time is defined as the time difference between when an image enters the queue and when it is first reviewed by a clinician. In a real-world clinic, it is typically the time from patient scan completion to the first case-open time by a radiologist. Wait-time-savings and delay are then defined as the wait-time difference between the without-AI workflow (control arm) and the with-AI workflow (test arm). 

\subsection{Workflows}
\begin{figure}[H]
\centering
\includegraphics[width=1\textwidth]{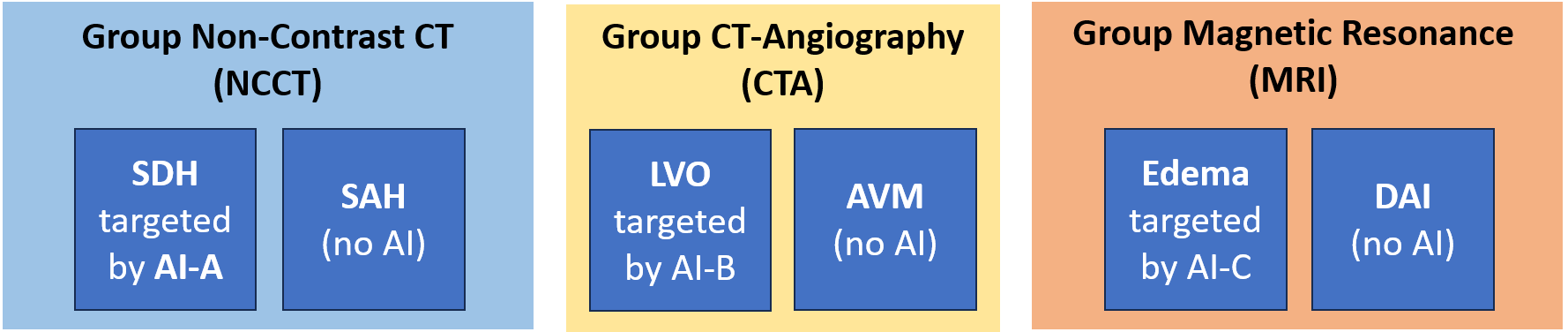}
\caption{\label{fig:Groups_Diseases_AIs}
An example of how a multi-AI, multi-disease workflow is defined. Bigger boxes represent different image groups by imaging modalities, whereas smaller blue boxes are a set of disease conditions detected using the corresponding imaging modalities. A set of AIs (AI-A to AI-D) are included to identify and triage a subset of disease conditions. NCCT = Non-Contrast Computed Tomography; CTA = Computed Tomography Angiography; MRI = Magnetic Resonance Imaging; SDH = Subdural Hemorrhage; SAH = Subarachnoid Hemorrhage; LVO = Large Vessel Occlusion; DAI = Diffuse Axonal Injury; AVM = Arteriovenous malformations; AI = Artificial Intelligence.}
\end{figure}

In this section, we introduce the workflows considered in this work and define the relevant parameters.
In general, a workflow consists of a set of image groups $M$, a set of disease conditions $B$, a set of AI-triage devices $A$, and a set of priority classes $C$. Figure \ref{fig:Groups_Diseases_AIs} is an example of how the three sets define a workflow.
\begin{itemize}
    \item Let $L$ be the total number of image groups i.e. $M = \{m_1, m_2, \cdots, m_L\}$. An image group $m_l$ is defined as a group of images that satisfy the inclusion criteria of an AI-triage device or a group of images that no AI-triage device reviews. For example, because an AI-triage typically checks for modality (e.g. via image header) or anatomy (e.g. head, chest, etc.) before processing an image, magnetic resonance (MR) images may be considered as one image group, whereas Computed Tomography (CT) images would be considered another group that may be further separated into with-contrast CT-Angiography (CTA) or non-contrast CT (NCCT).
    \item A workflow may include patients of $N$ different disease conditions: $B = \{b_1, b_2, \cdots, b_N\}$. Disease conditions can be ordered so that $b_1$ is most time-sensitive and $b_N$ is least time-sensitive. Each disease must belong to exactly one image group, but an image group may have patient images with multiple diseases. For model simplicity, a patient can have \emph{at most one} disease condition in each image group.
    \item Each disease condition $b_i$ may have an AI-triage device $a_i$ for identification and prioritization:  $A= \{a_1, a_2, \cdots, a_N\}$. An AI is defined by its diagnostic performance (namely sensitivity and specificity) to identify the corresponding disease condition. In this work, only the device performance is simulated, and no actual AI development was done. Hence, independent AI performance among conditions is assumed. While every AI-triage corresponds to a specific disease condition, not all conditions have an AI-triage. For example, in Fig. \ref{fig:Groups_Diseases_AIs}, no AI are incorporated to identify and triage cases with SAH. In this case, the corresponding $a_i$ is set to `None'. 
    \item Depending on the protocol, a workflow with $K$ AIs for prioritization might have $K+1$ priority classes: $C = \{c_1, c_2, \cdots, c_{k\leq N}, c^-\}$, where $c^-$ indicates the AI-negative class.  $c_k$ is ranked from the highest priority to the lowest, and the number of AI-positive priority classes is the equal to the number of $a_i$ that is not `None'. For example, in a workflow with five disease conditions and two AIs (one for condition 3 and another for condition 5), $C=\{c_1, c_2, c^-\}$ where $c_1$ are cases triaged by $a_3$, and $c_2$ are those flagged by $a_5$.
\end{itemize}

\begin{figure}[H]
\includegraphics[width=1\textwidth]{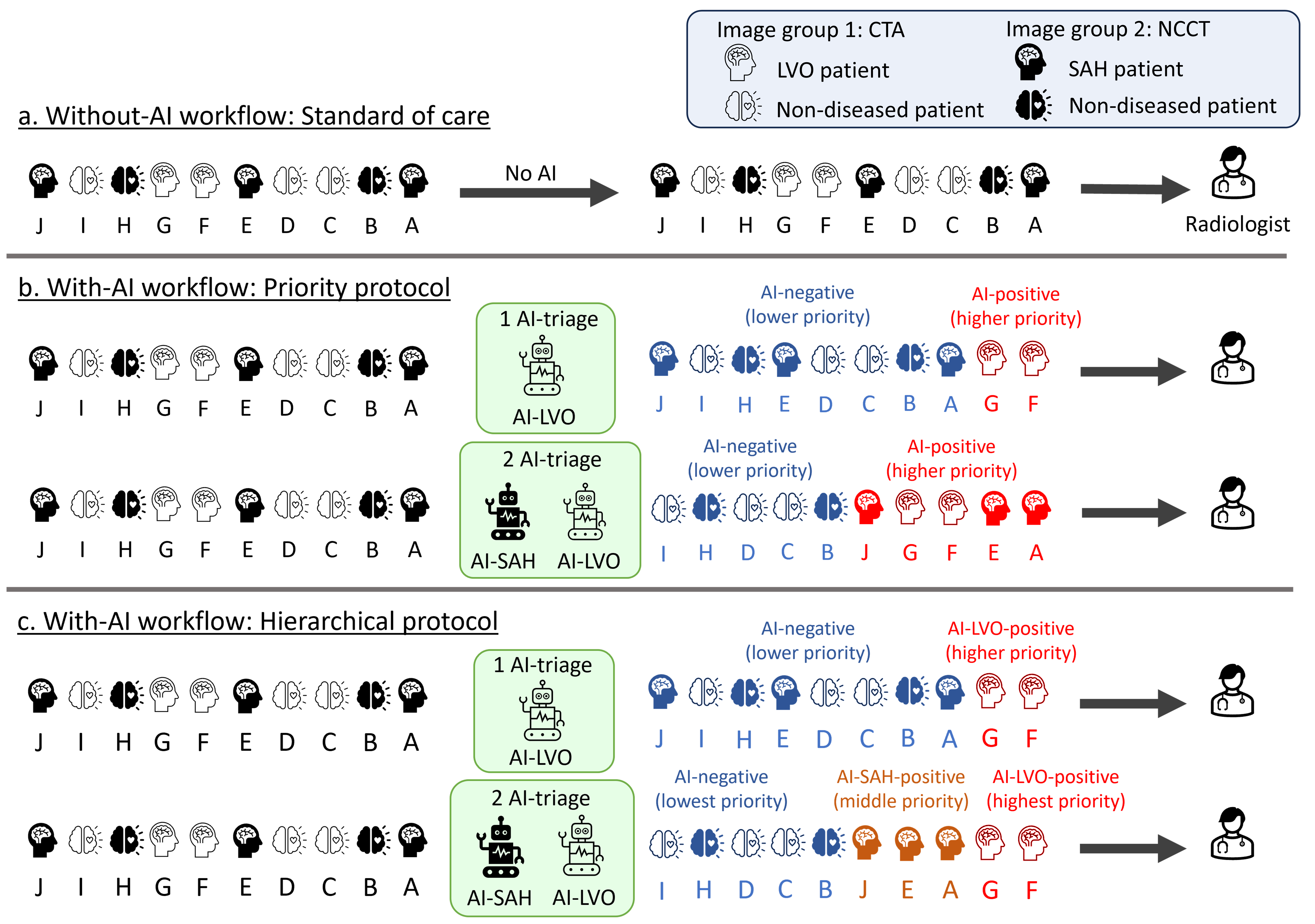}
\caption{\label{fig:workflowDefinition}
 An example of the types of multi-AI, multi-disease workflow. a) A without-AI workflow has a first-in first-out (FIFO) queue in which in-coming images are reviewed in the order of their arrivals. b) In a with-AI, priority workflow, all AI-positive cases are triaged as one high-priority class (red), and AI-negative cases are in the low-priority class (blud). Within the prioritized class, images are reviewed in the order of their arrivals regardless of the time-sensitiveness among disease conditions. c) A with-AI, hierarchical workflow further prioritizes the AI-positive cases by their time-sensitiveness based on a disease hierarchy (i.e. LVO being more time-sensitive than SAH) defined by the user. Hence, the 2-AI hierarchnical workflow has three priority classes: cases flagged by AI-LVO have the highest priority (red), cases flagged by AI-SAH have a middle priority (brown), and cases that are deemed negative by both AIs have the lowest priority (blue). SAH = Subarachnoid Hemorrhage; LVO = Large Vessel Occlusion; AI = Artificial Intelligence.}
\end{figure}

Figure~\ref{fig:workflowDefinition} shows the types of workflow considered in this work using patient brain images as an example. A radiologist reading queue consists of patients with two disease conditions $B=\{\text{LVO}, \text{SAH}\}$, as well as healthy patients without any brain diseases. LVO is typically diagnosed using CT-angiography (CTA) imaging, whereas diagnoses of SAH is usually based on non-contrast CT (NCCT) images. Hence, this example has only two image groups $M= \{\text{CTA}, \text{NCCT}\}$. 

When no AI-triage device is used, patient images are reviewed in the order of their arrival time i.e. first-in first-out (FIFO) (see Fig\ref{fig:workflowDefinition}a). The without-AI workflow is considered as the standard-of-care and the control arm in this work.

There is no universal standard on how multiple AI-triage devices are integrated into a radiologist workflow. Therefore, we consider two possible protocols:

\begin{itemize}
\item In a \emph{priority} protocol, all AI-positive patients from any AI $a_i$ are prioritized in the front of the queue, without distinguishing between specific AIs or disease conditions in the prioritization. Thus, a with-AI, priority workflow consists of two priority classes: AI-positive and AI-negative patient images. Figure\ref{fig:workflowDefinition}b illustrates the distinction between having one AI-triage device and having two, in which AI-positive cases are read in a FIFO ordering.
\item In a \emph{hierarchical} protocol, disease conditions are ordered by their time-sensitiveness, and patient images labeled as positive by AIs targeting higher-priority conditions are read first. Thus, a hierarchical workflow has at most $N+1$ priority classes: $a_1$-positive, $a_2$-positive, $\cdots, a_n$-positive, and all AI-negative patient images. If a clinic considers LVO to be more time-critical than ICH, a 2-AI hierarchical workflow prioritizes AI-identified LVO patient images in front of ICH patient images (as shown in Figure\ref{fig:workflowDefinition}c).
\end{itemize}

For each prioritization protocol, two queueing disciplines are considered: preemptive-resume and non-preemptive. In a preemptive-resume discipline, a radiologist reading a lower-priority case immediately stops the case in-hand and jumps to the higher-priority case upon its arrival.In a non-preemptive discipline, the radiologist finishes the lower-priority case first before starting the review of the in-coming higher-priority case. 

For model simplicity, this study focuses on workflows involving a single radiologist reviewing images. While the simulation can simulate workflows with multiple radiologists, the 1-radiologist workflow allows us to estimate the maximum potential time savings and delays for various subgroups.

\subsection{Parameters}
To model time-savings and delay for patients in these workflows, the following parameters are defined.
\begin{itemize}
\item Group probabilities $g_i$ are defined as the fraction of patient images in group $m_i$ among all patient images.
\item $\lambda$ is the Poisson arrival rate of all patient images. Patient images can be divided into subgroups, and each subgroup $i$ has Poisson arrival rate $\lambda_i = q_i \lambda$, where $q_i$ is the fraction of image subgroup $i$ with respect to all patient images. For example, the arrival rate of cases in an image group depends on its group probability. If we are interested in the arrival rate of patient images with condition A prioritized by any AIs, $q_i$ will depend on the prevalence of condition A as well as the sensitivity and specificity of the AIs involved in the workflow.
\item Disease prevalence rates are defined as $\pi_1, \pi_2, \cdots, \pi_N$, where $\pi_i$ corresponds with disease $b_i$:
$$\pi_i = \frac{\text{Number of cases with disease } b_i}{\text{Total number of cases in the queue}}$$
\item For an AI $a_i$ targeting at disease $b_i$, its sensitivity $\mathrm{Se}_i$ and specificity $\mathrm{Sp}_i$ are calculated with respect to number of cases processed by the AI. By definition of image group, the denominator is the number of patient images within the image group that the AI belongs to. 
$$\mathrm{Se}_i = \frac{\text{Number of AI-positive cases processed by AI } a_i \text{ with disease } b_i}{\text{Number of cases processed by AI } a_i \text{ with disease } b_i},$$
and
$$\mathrm{Sp}_i = \frac{\text{Number of AI-negative cases processed by AI } a_i \text{ without disease } b_i}{\text{Number of cases processed by AI } a_i \text{ without disease } b_i}.$$
When there is no AI targeting at disease $b_i$, the effect to the case prioritization is equivalent to having an AI $a_i$ that labels every image in the image group as $a_i$-negative for $b_i$. Hence, if $a_i = \mathrm{None}$, then by default $\mathrm{Se}_i = 0$ and $\mathrm{Sp}_i = 1$.
\item The radiologist mean reading times are denoted by $S_i$'s. In queueing theory, the exponential distribution is the simplest model for reading time, where reading rates $\mu_i$ are the inverse of mean reading times i.e. $\mu_i = 1/S_i$. The average time for reviewing a patient image with disease $i$ is denoted by $S_{i, \D}$, while $S_{i, \N}$ is the mean reading time for a patient image that does not have disease $i$. If diseases $b_i$ and $b_j$ are in the same image group, then $S_{i, \N} = S_{j, \N}$. 
\item $\rho$ is the traffic intensity defined as $\rho = \lambda/\mu_{\text{eff}}$, where $\mu_{\text{eff}}$ is the effective reading rate considering all priority classes. $\rho$ ranges from 0 (with no arriving images) to 1 (a very congested clinic).
\item $p_i^{+}$ is the probability that (with respect to all images in the queue) an image belongs to the $a_i$-positive group but is negative for all higher-priority AIs in the same group (else, the image will be added to the higher-priority group). Similarly, $p_i^{-}$ is the probability that an image belongs to the set of AI-negative patients in the image group that $a_i$ belongs to.

\end{itemize}

\section{Methods}
We applied both theoretical calculation and simulation approach to predict mean time-savings and delays for different subgroup of patient images in the workflows of various settings. Mean time-savings and delays are based on the wait-time difference between with-AI and without-AI workflows. By definition, a positive difference for a subgroup of patient images implies that they are on average delayed due to the introduction of AI-triage devices (risks). If it is negative, the patient images in the subgroup experience a net time-savings (benefits). 

\subsection{Theoretical Calculations}
\label{sec:theory}

This section focuses on the wait-time calculations in various with-AI workflows. The theoretical framework for calculating mean wait-times in a without-AI workflow was outlined in \cite{Thompson}. Although \cite{Thompson} only considers single-AI, single-disease workflows, the \textit{mean} wait-time across \textit{all} patients in a without-AI workflow is the same for multi-AI, multi-disease workflows due to the randomness in Poisson arrival pattern under the assumption that AI does not affect reading times.

Wait-time calculations for AI-positive and AI-negative subgroups in the four with-AI workflows are presented. Sections \ref{subsec:Ppriority} and \ref{subsec:Phierarchical} extend the work of \cite{Thompson} to analyze wait-times for various subgroups under the preemptive discipline given priority and hierarchical protocols. Methods to obtain theoretical estimates for wait-times in the non-preemptive setting are derived from \cite{Raicu} and are described in Sections \ref{subsec:NPpriority} and \ref{subsec:NPhierarchical} for priority and hierarchical protocols respectively. The derivation of the mean wait-times for each diseased and non-diseased subgroup are obtained directly from the AI-positive and AI-negative estimates; details of this conversion are left to Appendix~\ref{App:waittime}.

\subsubsection{Preemptive Priority Workflow}
\label{subsec:Ppriority}
In a preemptive workflow, a radiologist reviewing a lower-priority patient image is immediately interrupted by a higher-priority case upon its arrival. Since a priority protocol has only two priority classes (AI-positive and AI-negative), Model 2 in \cite{Thompson} can be re-used for wait-time calculations in a with-AI workflow under the priority protocol and preemptive discipline setting. The workflow in Model 2 is limited to only one radiologist but allows different reading rates for the diseased and non-diseased patient image subgroups. In a 1-AI, 1-disease scenario, Model 2 takes as inputs the mean reading rates ($\mu_{+}$,  $\mu_{-}$) and the arrival rates ($\lambda_{+}$, $\lambda_{-}$) of the AI-positive and AI-negative classes. These rates depend on probabilities such as disease prevalence and AI sensitivity and specificity.

When a workflow involves multiple AI-triage devices, and all AI-positive cases are triaged as one high-priority class, the queue can be viewed as a 2-priority-class system with effective rates ($\mu_{+}$,  $\mu_{-}$, $\lambda_{+}$, $\lambda_{-}$) dependent on a set of prevalence rates, AI sensitivities and specificities, and group probabilities. For example, since the AI-positive subgroup consists of triaged cases with different conditions $b_i$ and non-diseased cases, and each diseased and non-diseased subgroup has its own Poisson arrival rates $\lambda_i$ (i.e. the inter-arrival time between two consecutive cases in the same subgroup $i$ is exponentially distributed), the inter-arrival time between two AI-positive cases are hyper-exponentially distributed. The effective rates $\lambda_{\pm}$ are therefore
\begin{equation}
\frac{1}{\lambda_{\pm}} = \sum_{i=1}^K \sum_{\substack{ j \in B_i}} \frac{\mathbb{P}(\text{diseased } b_j \vert a_i^{\pm} \text{ subgroup})}{\lambda_i} + \frac{\mathbb{P}(\text{non-diseased}\vert a_i^{\pm} \text{ subgroup})}{\lambda_\text{ND}},
\end{equation}
where $B_i$ is a set of disease conditions in the image group that $a_i$ belongs to. This expression involves conditional probabilities, such as the probability that a patient image has disease $b_i$ given that it belongs to the AI-positive class; these calculations are left to Appendix~\ref{App:condprob}. Once the effective rates are determined, the calculations for mean wait-times in the AI-positive and AI-negative groups follow from the methods described in \cite{Thompson}.

\subsubsection{Preemptive Hierarchical Workflow}
\label{subsec:Phierarchical}
Unlike the priority protocol, in a hierarchical protocol, disease conditions are ordered by time-sensitiveness with a lower rank index being more time-critical. For example, cases triaged by $a_1$ has a higher priority than those by $a_2$, etc. Under a preemptive discipline, a radiologist reading an $a_i$-positive image may be interrupted only by images deemed AI-positive by $a_{j < i}$ and will not be interrupted by AI-positive cases from $a_{j\geq i}$. Readings of AI-negative images may be interrupted by AI-positive images from any AIs in the workflow. 

Hence, the total wait-time for $a_i$-positive cases depends only on the cumulative wait-time of AI-positive images triaged by $a_{j<i}$. This set of cases (denoted as $H$) consists of AI-positive images from $a_1, a_2, \cdots, a_{i-1}$. Similarly, the total wait-time for $a_{i+1}$-positive cases depends on the cumulative wait-time of images triaged by $a_1, a_2, \cdots, a_i$ (denoted as $L$). The mean total wait-time $\overline{W_{a_i}^{+}}$ for $a_i$-positive cases can then be written as the difference in total wait-time between the two sets of AI-positive patients:
$$\overline{W_{a_i}^{+}} = \overline{W^{+}_{a_j=1,\cdots,i}}-\overline{W^{+}_{a_j = 1, \cdots, i-1}}=\overline{W^{+}_L}-\overline{W^{+}_H}.$$

Computation of the mean wait-time ($\overline{W^+_L}$,  $\overline{W^+_H}$) are challenging for a general workflow with different reading rates among diseased and non-diseased subgroups. However, by limiting our calculation to equal mean reading rates for all subgroups, we can re-use Model 1 in \cite{Thompson} to calculate mean wait-times for the two sets of cases ($L$ and $H$) based on the probabilities $\pi^+_{\{L,H\}}$ that an image is AI-positive from one of the AIs in the set. The calculation of $\pi^+_{\{L,H\}}$ is provided in Appendix~\ref{App:probPi}. By weighting the mean wait-times by the prevalence rates and taking their difference, the average wait-time for $a_i$-positive images $\overline{W_{a_i}^{+}}$ can be expressed as
$$\overline{W_{a_i}^{+}} = (\overline{W_L^{+}} \pi_{L}^{+} - \overline{W_{H}^{+}} \pi_{H}^{+})/\pi_{a_i}^{+}.$$

The calculation of mean per-patient wait-time $\overline{W^{-}}$ for AI-negative images is the same as the one in Model 1 in \cite{Thompson}  irrespective of the multiple AIs. This can be computed by combining all AIs and the corresponding image and disease groups to obtain their ``effective'' combined AI performance and prevalence, which are used as inputs to the theoretical calculations for the AI-negative class in Model 1 in \cite{Thompson}. 

The above calculation is repeated $N$ times from the most time-critical disease condition to the least time-sensitive condition. At each step, the $L$ and $H$ sets are redefined, and the probabilities $\pi^+_{\{L,H\}}$ are recalculated. The average wait-times for all diseases are therefore computed iteratively using the above equation. 

\subsubsection{Non-Preemptive Priority Workflow}
\label{subsec:NPpriority}

In a non-preemptive discipline, a radiologist reviewing a lower-priority patient image is \emph{not} immediately interrupted by a higher-priority case upon its arrival. Instead, the radiologist completes the lower-priority case in-hand first before reviewing the new higher-priority case.

With only two priority classes under the priority protocol (AI-positive from any AIs and AI-negative from all AIs), the mean wait-times for the AI-positive class can be calculated based on a Eq. 2 from \cite{Raicu}.
\begin{equation}\label{eq:NPpriorityPos}
\overline{W^{+}} = \frac{1}{\mu_+}\frac{\rho_+}{1-\rho_+}.
\end{equation}
$\lambda_+$ and $\rho_+$ are the effective arrival rate and traffic intensity for the overall AI-positive cases flagged by any AIs. Similarly, Eq. 7 from \cite{Raicu} can be used to calculate the mean wait-time for AI-negative class.
\begin{equation}\label{eq:NPpriorityNeg}
\overline{W^-} = \frac{1}{\mu_-} \frac{\frac{\mu_-}{\mu_+} \cdot \frac{\rho_+}{1-\rho_+} + \rho}{1-\rho},
\end{equation}
where $\rho=\rho_+ + \rho_-$. The derivations of the conditional probabilities that are involved in calculating the effective arrival rate and traffic intensity are left in Appendix~\ref{App:condprob}.

\subsubsection{Non-Preemptive Hierarchical Workflow}
\label{subsec:NPhierarchical}

In a hierarchical protocol, disease conditions are ordered by time-sensitiveness. Unlike the priority protocol with only two priority classes, the hierarchical protocol has a set of priority classes $C$.

Equation 21 in \cite{Raicu} is applicable to calculating mean wait-time for each priority class in a multi-priority-class workflow. Specifically, the mean wait-time of the $k$-th priority class ($\overline{W_k^{+}}$) depends on the arrival rates $\lambda$, mean reading times $S$, and the variances $V$ of the reading time distributions from all $K$ AI-positive priority classes~\cite{Raicu}:
\begin{equation}\label{eq:NPhierarchicalPos}
\overline{W_k^{+}} = \frac{\sum_{i=1}^K \lambda_i V_i + \lambda_- V_-}{2 \Big(1-\sum_{j=1}^{k-1} \lambda_j S_j\Big) \Big(1-\sum_{j=1}^k \lambda_j S_j\Big)}.
\end{equation}
Similarly, the mean wait-time in for the AI-negative class is 
\begin{equation}\label{eq:NPhierarchicalNeg}
\overline{W^-} = \frac{\sum_{i=1}^K \lambda_i V_i + \lambda_- V_-}{2 \Big(1-\sum_{j=1}^{K} \lambda_j S_j\Big) \Big(1-\sum_{j=1}^K \lambda_j S_j - \lambda_{-} S_{-}\Big)}.
\end{equation}

The arrival rate of the $i$-th priority class depends on the probabilities $p_i^+$ that an image is $a_i$-positive but negative for all higher-priority AIs in the same group. The calculation can be found in the Appendix~\ref{App:probPi}. With $p_i^+$, the arrival rate of images in $i$-th priority class is
$$\lambda_i = p_i^{+} \lambda.$$
Similarly, for AI-negative images in the lowest priority class,
$$\lambda_{-} = p_{-} \lambda = \big(1-\sum_{k=1}^K p^+_k\big) \lambda.$$

Besides arrival rates, Eq.~\ref{eq:NPhierarchicalPos} and \ref{eq:NPhierarchicalNeg} also need the means and variances of the reading time distributions for all $K$ priority classes. Because the $k$-th priority class can consist of different diseases, and the reading times of diseased and non-diseased subgroups are exponentially distributed with different means, the reading time distribution for the $k$-th priority class is hyper-exponential. Let $B_i$ be the set of disease conditions in the image group containing AI $a_i$; for example, in Figure \ref{fig:Groups_Diseases_AIs}, $B_1$ corresponding to the group that AI-A belongs to includes SDH and SAH. The mean of the hyper-exponential reading time distribution for the $k$-th priority class is then
\begin{align*}
S_k = \sum_{b_j \in B_i} &\mathbb{P}(\text{diseased } b_j \vert a_i^{+} \text{ subgroup}) \cdot S_{j, \D}\\
& + \mathbb{P}(\text{non-diseased} \vert a_i^{+} \text{ subgroup}) \cdot S_{j, \N},
\end{align*}
where $S_{j, \D}$ and $S_{j, \text{ND}}$ are the mean read-time for cases with condition $j$ and that for non-diseased cases in the same image group. The second moment $V_k$ of the distribution is
\begin{align*}
V_k = \sum_{b_j \in B_i} &\mathbb{P}(\text{diseased } b_j \vert a_i^{+} \text{ subgroup}) \cdot 2\cdot S_{j, D}^2\\
& + \mathbb{P}(\text{non-diseased} \vert a_i^{+} \text{ subgroup}) \cdot 2 \cdot S_{j, \N}^2.
\end{align*}

The AI-negative group consists of images deemed negative by all AIs from all image groups $M$. Its mean reading time is thus
\begin{align*}
S_{-} = \sum_{m_i\in M} \sum_{b_j \in B_i} &\mathbb{P}(\text{diseased } b_j \vert a^{-} \text{ subgroup})  \cdot S_{j, \D} \nonumber \\
&+  \mathbb{P}(\text{non-diseased, group } m_i \vert a^{-} \text{ subgroup})  \cdot S_{i, \N}
\end{align*}
with second moment
\begin{align*}
V_{-} = \sum_{m_i\in M} \sum_{b_j \in B_i} &\mathbb{P}(\text{diseased } b_j \vert a^{-} \text{ subgroup}) \cdot 2\cdot S_{j, \D}^{2} \nonumber \\ &+ \mathbb{P}(\text{non-diseased, group } m_i \vert a^{-} \text{ subgroup}) \cdot 2\cdot S_{i, \N}^2.\\
\end{align*}

\subsection{Simulation}
\label{sec:simulation}
The simulation software described in \cite{Thompson} is extended to accommodate multi-AI, multi-disease workflows. To dynamically organize the different AI/disease configurations, the simulation software includes a disease tree and a disease hierarchy to accept any number of image groups, disease conditions, and AI-triage devices. Each image group is defined by its name, its probability with respect to all cases in the workflow, as well as the disease condition(s) that belong to this group. Each disease condition is defined by its name, its rank (for hierarchical protocol), as well as the prevalence with respect to the cases in the group. Each AI is defined by its name, the image group it belongs to, the target disease, and its sensitivity and specificity. The simulation of patient image flow is modified accordingly to incorporate the non-preemptive discipline and hierarchical protocol.

Each simulated patient image arrives according to a Poisson distribution and is randomly assigned characteristics including the image group it belongs to, the disease condition, the disease status, and AI-call status from all the AIs within the image group (if any). Each image is then simultaneously placed in two parallel worlds: in the without-AI world, images in the reading queue are ordered via FIFO, whereas images in the with-AI world are put in different priority levels depending on the user-specified protocol. If it is a priority protocol, each image in the with-AI world is either in high priority (AI-positive by any AIs) or low priority (AI-negative by all AIs). If it is a hierarchical protocol, the priority level of the image is based on the AI-positive label of the most time-critical disease condition. For example, consider a hierarchical queue with two AIs triaging two (A and B) of the three (A, B, and C) conditions where A > B > C in their time-sensitiveness ranking. If a patient image is diseased with condition B and labeled as positive by both AIs, its priority level $k$ is 1 (highest-priority). If the image is labeled as AI-positive only by AI-B, then $k$=2. If none of the AIs flag it, then the image belongs to the AI-negative class with the lowest priority.

The simulation tracks waiting times and queue compositions as well as calculating the difference in waiting times between the two worlds to evaluate whether the AIs provides time savings or delays for each patient subgroups. While the two arms are either with all AIs or without any AIs, if one is interested in the time-savings/delay due to one subject AI-triage device, one can run the software twice (one with the subject AI and one without) to estimate the time difference due to the subject AI-triage device.

It should be noted that the simulation software is designed to estimate mean wait-time savings and delays for different subgroups in a much more complex workflow. Due to the limitation of the theoretical approach, calculations in \S~\ref{sec:theory} are restricted to one radiologist reading images and, in some cases, equal reading rates between diseased and non-diseased subgroups. Simulation does not have those limitations and is capable of simulating radiologist workflows in the real-world where these assumptions may not be applicable. The full software suit can be found in a publicly available GitHub repository~\cite{ThompsonMultiQuCAD}.

\section{Experimental Scenarios}
\label{sec:workflows}
To understand how deploying multiple AI-triage devices affect patient wait-time, we defined three experimental workflow scenarios, all using brain CTA and NCCT images as a representative example. All three experimental scenarios involve images of patients with LVO, SDH, and/or SAH, along with healthy brains, but they differ by the AI-triage device(s) involved in re-arranging the cases. Our goal is to compare the trade-offs between time-savings and delay for images with high- versus low-priority disease conditions in these three brain-related experimental scenarios. Agreement between theoretical calculation and simulation is also reported as a function of wait-time difference between each workflow scenario and the without-AI workflow.

To further assess the consistency between simulation and theory in more complex workflows, we hypothesized an additional scenario with nine disease conditions and four AI-triage devices and compared the theoretical calculation and simulation results.

\subsection{Workflows}
\label{subsec:workflows} 

\begin{figure}[H]
\includegraphics[width=1\textwidth]{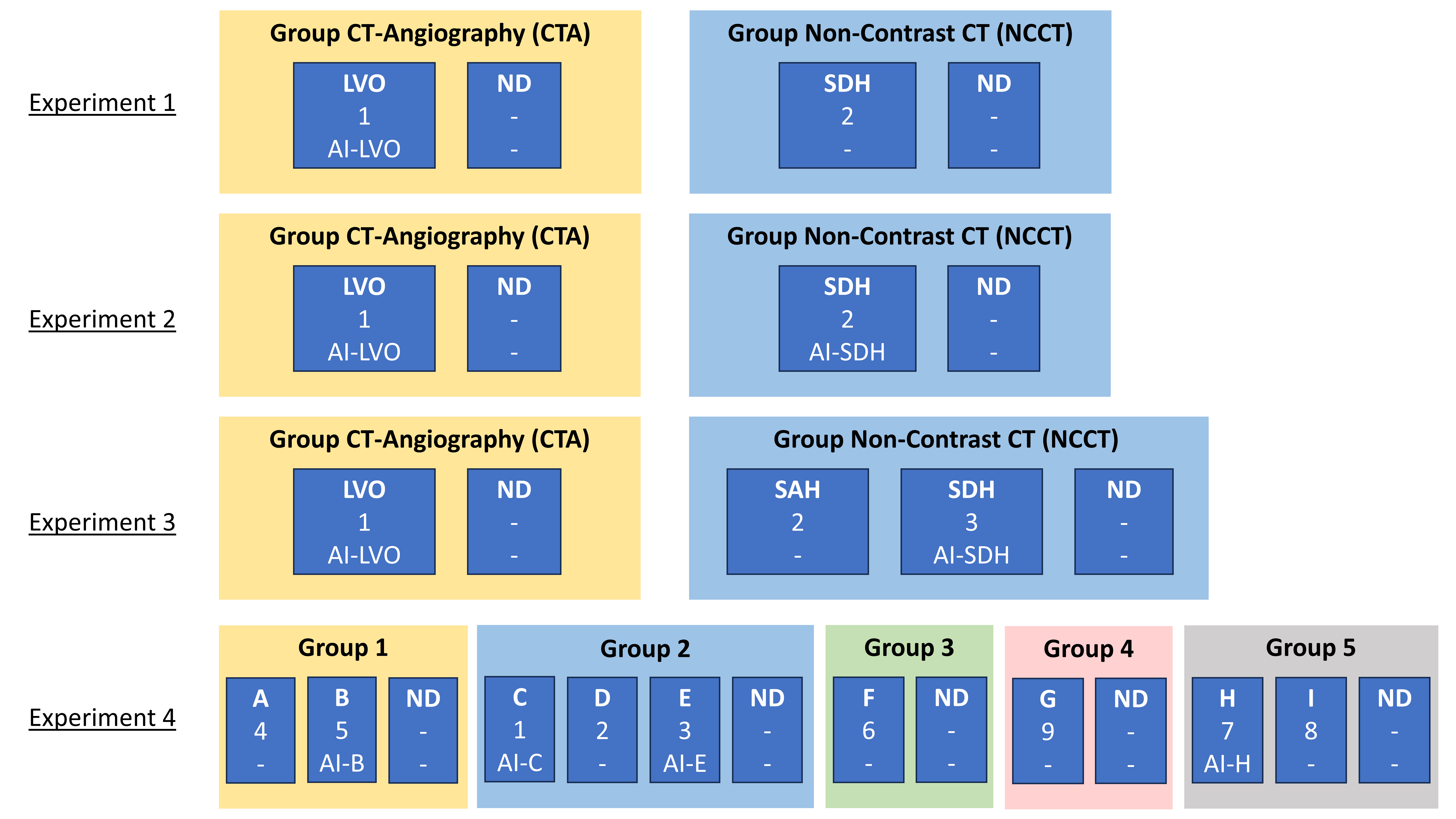}
\caption{\label{fig:expWorkflowDefinition}
The AI and disease condition settings for the four experimental workflow scenarios. Each box with different shade colors represents an image group. Each group consists of patients with a disease condition as well as non-diseased (ND) patients. The smaller blue boxes represent the disease condition and ND subgroups. Each disease condition is defined by a disease name, its time-sensitiveness ranking (with 1 being the most time-critical condition), and the name of the corresponding AI-triage device (if any). LVO = Large Vessel Occlusion; SAH = Subarachnoid Hemorrhage; SDH = Subdural Hemorrhage; AI = Artificial Intelligence.}
\end{figure}

Figure~\ref{fig:expWorkflowDefinition} illustrates the image groups and disease conditions in the four experimental workflow scenarios. Each of them has a unique multi-AI, multi-disease workflow. This section defines the disease conditions and AI-triage devices in each workflow scenario.

\subsubsection{Experiment 1: LVO and SDH with one AI}
\label{ExpScenarios:subsubsec:exp1}
We start with the simplest scenario workflow with images of patients with LVO or SDH, along patients with healthy (non-diseased) brains; SAH is not considered. As a demonstration, we hypothetically assign LVO to be more time-critical than SDH. Clinically, LVO is typically diagnosed via CTA~\cite{Czap}, and SDH is diagnosed via NCCT~\cite{Penckofer}. Therefore, this workflow scenario has two image groups. Only one AI-triage device (AI-LVO) is involved to identify patients with LVO and triage the AI-LVO-positive patients at the top of the reading queue. No AI is involved to identify SDH cases. Experiment 1 is a 1-AI, 2-disease workflow which can be used as a reference arm to study any additional patient wait-time-savings and delays when compared to Experiments 2 and 3.

\subsubsection{Experiment 2: LVO and SDH with two AIs}
\label{ExpScenarios:subsubsec:exp2}
The second scenario workflow also focuses on LVO and SDH but with an additional AI-triage device (AI-SDH) for identifying and prioritizing SDH cases. Since the two conditions belong to different image groups (LVO in CTA and SDH in NCCT), the two AIs process different sets of images in the workflow. Comparing wait-time-savings for LVO patient images (highest-priority in a hierarchical protocol) between Experiments 1 and 2 is of interest, as is examining the time-savings benefit to SDH cases when an AI-SDH device is introduced.

\subsubsection{Experiment 3: LVO, SDH, and SAH with two AIs}
Experiment 3 extends the complexity by including all three conditions: LVO, SDH, and SAH. Both SDH and SAH belong to the NCCT image group~\cite{Penckofer}, and we hypothetically assign a time-sensitive ranking such that LVO > SAH > SDH. LVO and SDH continue to benefit from AI-triage devices, but no AI is involved to identify patients with SAH. Without an AI for identifying and triaging SAH cases, most of these cases will likely have the lowest priority and experience time delay. However, since AI-SDH analyzes all NCCT images, it may mis-label and prioritize some SAH-present images depending on its False-Positive rate. With the addition of SAH, this workflow scenario aims to elucidate the complex and interrelated trade-offs between time-saving benefits and delay risks across diverse patient imaging cases.

\subsubsection{Experiment 4: a 4-AI, 9-disease workflow scenario}
Last, we considered a hypothetical workflow with four AI-triage devices trained to identify four of the nine disease conditions spread across five image groups. No specific clinical context is assumed. The goal here is to show agreement between theoretical calculation and simulation and to demonstrate generalization of this work to more complex workflow scenarios in any radiology department. 

\subsection{Parameters}
This section describes how the numerical values of parameters related to image group, disease conditions, AI performance, and clinical settings are determined. We first focused on the parameters in the three brain-related scenarios followed by those in the hypothetical 4-AI, 9-disease workflow scenario.

Parameters related to image groups and disease conditions (listed below) are fixed in this study to emulate the three experimental scenarios.  
\begin{itemize}
    \item For each group (CTA and NCCT)
    \begin{itemize}
        \item Group probability with respect to all cases
        \item Mean reading time for images with non-diseased brains
    \end{itemize}
    \item For each condition,
    \begin{itemize}
        \item Prevalence within the group
        \item Mean reading time 
        \item Time-sensitiveness ranking for hierarchical protocol
    \end{itemize}
\end{itemize}
To our best ability, we could not find any literature that studies the proportion of CTA and NCCT in a clinic. Therefore, CTA and NCCT group probabilities are set to 0.3 and 0.7 respectively, reflecting the fact that NCCT is performed more widely and frequently than CTA~\cite{Chung}. 

For disease condition-related parameters, we identified two publications that reported mean reading times and prevalence rates for LVO~\cite{Soun} and ICH~\cite{Savage}. The numerical values and references are shown in Table \ref{tab1}. The prevalence rates for the three conditions in Table~\ref{tab1} are used in our experimental scenarios. Since all the mean reading times are less than  5 minutes away from the total average, we set the read-time to be 30 minutes across all three conditions such that we can study the agreement between theoretical calculation and simulation across all workflow types including the preemptive, hierarchical workflow. Last, for demonstration purposes, we hypothetically ranked LVO as the most time-sensitive condition, followed by SAH, and then SDH as the least time-critical condition.
\begin{table}[H] 
\caption{Numerical values of disease-related parameters in the three experimental scenarios. \label{tab1}}
\begin{tabularx}{1\textwidth}{c|c|c|c}
\toprule
& \textbf{LVO}	& \textbf{SAH} & \textbf{SDH}\\
\midrule
\textbf{Prevalence}	& $\frac{48+47}{439+321} = 0.125$ {\fontsize{9}{10.8}\selectfont(p.3 \cite{Soun})} & 0.053 {\fontsize{9}{10.8}\selectfont(Table 1 \cite{Savage})} & $\frac{735+1368}{3716+628} = 0.21$ {\fontsize{9}{10.8}\selectfont(Table 1 \cite{Savage})} \\
\textbf{Read-time [min]} & 32.67 {\fontsize{9}{10.8}\selectfont(Table 1 \cite{Soun})} & 24.3 {\fontsize{9}{10.8}\selectfont(Table 5 \cite{Savage})}  & 24.3 {\fontsize{9}{10.8}\selectfont(Table 5 \cite{Savage})}\\
\bottomrule
\end{tabularx}
\end{table}



\begin{figure}[h!]
  \centering
    \includegraphics[width=1\textwidth]{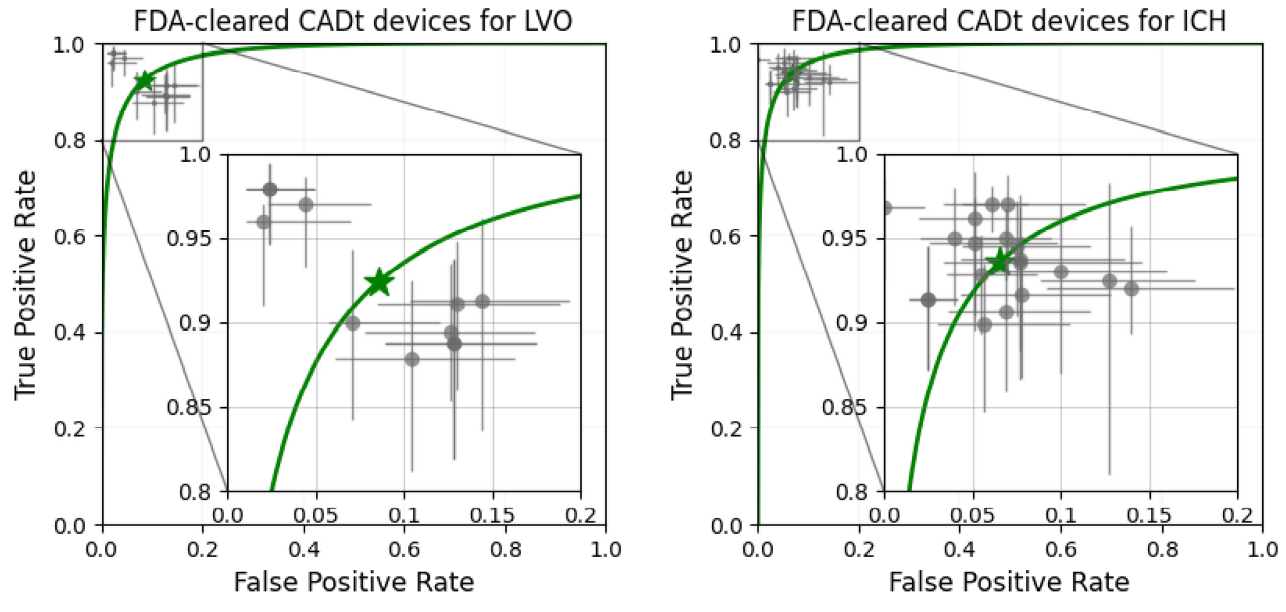}
  \caption{Sensitivity and specificity performance of FDA-cleared AI-triage devices targeting LVO (left) and ICH (right) respectively. Each gray dot is the sensitivity and specificity pair from an AI-triage device with their 95\% confidence intervals. The green star is the average performance among all devices. An ROC curve in green is obtained by fitting the average performance with a binomial assumption. ICH = Intracerebral Hemorrhage; LVO = Large Vessel Occlusion; FDA = U.S. Food and Drug Administration; ROC = Receiver Operating Characteristic; AI = Artificial Intelligence. \label{fig:FDA_CADt}}
\end{figure}

Which AI-triage device(s) are involved depends on the three experimental scenarios, and their default sensitivity (Se) and specificity (Sp) are set to the mean performance values among FDA-cleared AI-triage devices. We identified 12 AIs that target LVO (AI-LVO) and 12 targeting ICH (AI-ICH) from the public-facing OpenFDA API~\cite{openfda}. The Se and Sp of each AI are extracted from the corresponding 510(k) summary. As shown in Fig.~\ref{fig:FDA_CADt} (left), the average Se and Sp for LVO identification are 0.9236 and 0.9143 respectively, which are used when AI-LVO is involved in a workflow scenario. As for AI-SDH, the summary documents do not include the performance for the individual ICH subtype. Therefore, we used the average performance of Se = 0.9362 and Sp = 0.9343 (Fig.~\ref{fig:FDA_CADt}, right) from all available AI-ICHs. It should be noted that, in addition to the Se-Sp pair, the multi-QuCAD software also accepts the full ROC curve of an AI-triage device. More information can be found in the user manual of the software~\cite{ThompsonMultiQuCAD}.

In the three brain-related experimental scenarios, we investigated wait-time impacts in various patient image subgroups due to two factors: traffic intensity $\rho$ and AI performance. To study effects due to traffic intensity, we varied $\rho$ from 0 (empty queue) to 1 (very congested queue), using the image group and disease related parameters and the default AI sensitivity and specificity values defined above. For the impacts due to AI performance, we constructed a Receiver Operating Characteristic (ROC) curve for each AI using the average (Se, Sp) points under a binormal assumption, as shown in Fig.~\ref{fig:FDA_CADt} . We then varied the AI operating point along the curve, assuming a traffic intensity $\rho$ of 0.8 (which represents a moderately busy clinic) while fixing the image group and disease related parameters.

\begin{figure}[h!]
  \centering
    \includegraphics[width=0.85\textwidth]{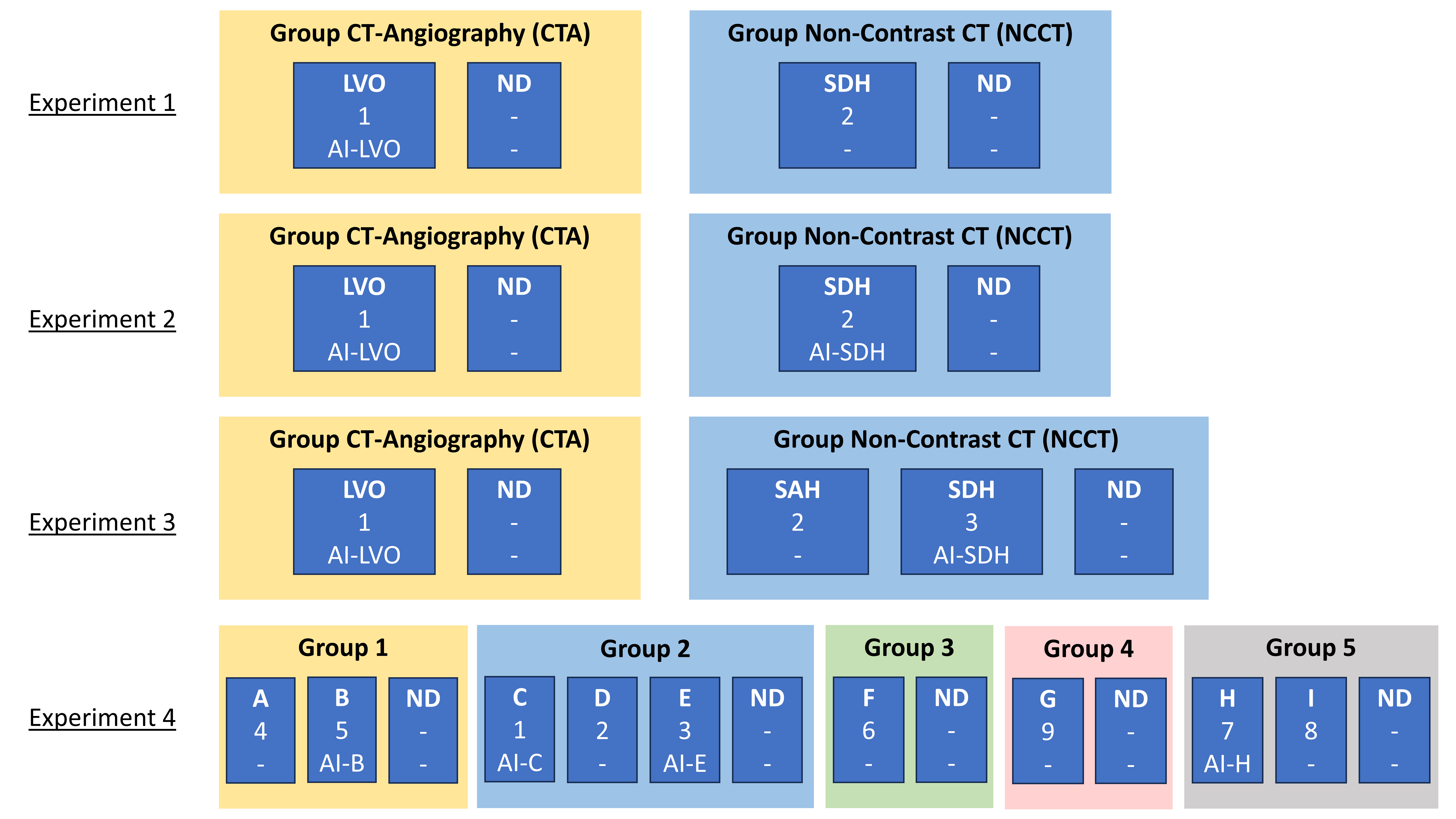}
  \caption{The numerical values of the parameters for the 4-AI, 9-disease workflow scenario (Experiment 4). (Left) Groups are defined by a block with its name, group probability, disease conditions involved, and the mean read-time of the non-diseased cases in the group. Each disease is defined by four parameters: disease name, its hierarchical rank, prevalence within the group, and the mean read-time. (Right) Each AI is defined by a block with its name, the group that the AI belongs to, the target disease condition, the true-positive fraction (TPF), and the false-positive fraction (FPF). Sensitivity is the same as TPF, whereas specificity is 1-FPF. A Receiver operating characteristic (ROC) curve can also be provided via a csv file.  \label{fig:hypoParams}}
\end{figure}

As for the hypothetical 4-AI, 9-disease workflow scenario in Experiment 4, the traffic intensity $\rho$ is fixed at 0.8. The values of other input parameters are shown in Fig.~\ref{fig:hypoParams} and are randomly chosen without assuming any clinical applications. All parameters are fixed, and only simulation-theory agreement is reported.

\subsection{Analysis}

For each of the experimental workflow scenarios and for each parameter setting, we calculated the theoretical mean wait-time differences in different disease condition subgroups for all four workflow settings: preemptive priority, preemptive hierarchical, non-preemptive priority, and non-preemptive hierarchical. In parallel, 1000 simulation trials were generated, each of which has 10,000 simulated patients. Each simulated patient image has a wait-time difference between the with-AI and without-AI workflow. Mean wait-time differences were stratified by patients with different disease conditions, and the 95\% confidence intervals (C.I.s) were calculated from the distributions of the mean wait-time difference.

For each of the three brain-related experimental workflow scenarios, we reported the wait-time difference (both theoretical and simulation approaches) as a function of traffic intensity and AI performance. The trade-offs between time-savings and delay among different patient subgroups are studied. 

For agreement study, we compared results from simulation and theoretical calculations via relative errors (RE):.
$$ \text{Relative Error} = \frac{\text{theory} - \text{simulation}}{\text{theory}},$$
where $\text{theory}$ and $\text{simulation}$ refer to wait-time difference from theoretical calculations and simulation results respectively. Only the RE for Experiment 3 are reported across traffic intensity, disease prevalence, and mean reading time because this workflow scenario is the most complicated among Experiments 1 to 3. For the 4-AI, 9-disease workflow, we directly reported the absolute mean wait-time and wait-time difference for different disease condition subgroups.

\section{Results}
\label{sec:results} 

\subsection{Experiment 1: LVO and SDH with one AI}

\begin{figure}[h!]
\begin{tabular}{ll}
    \includegraphics[width=0.47\textwidth]{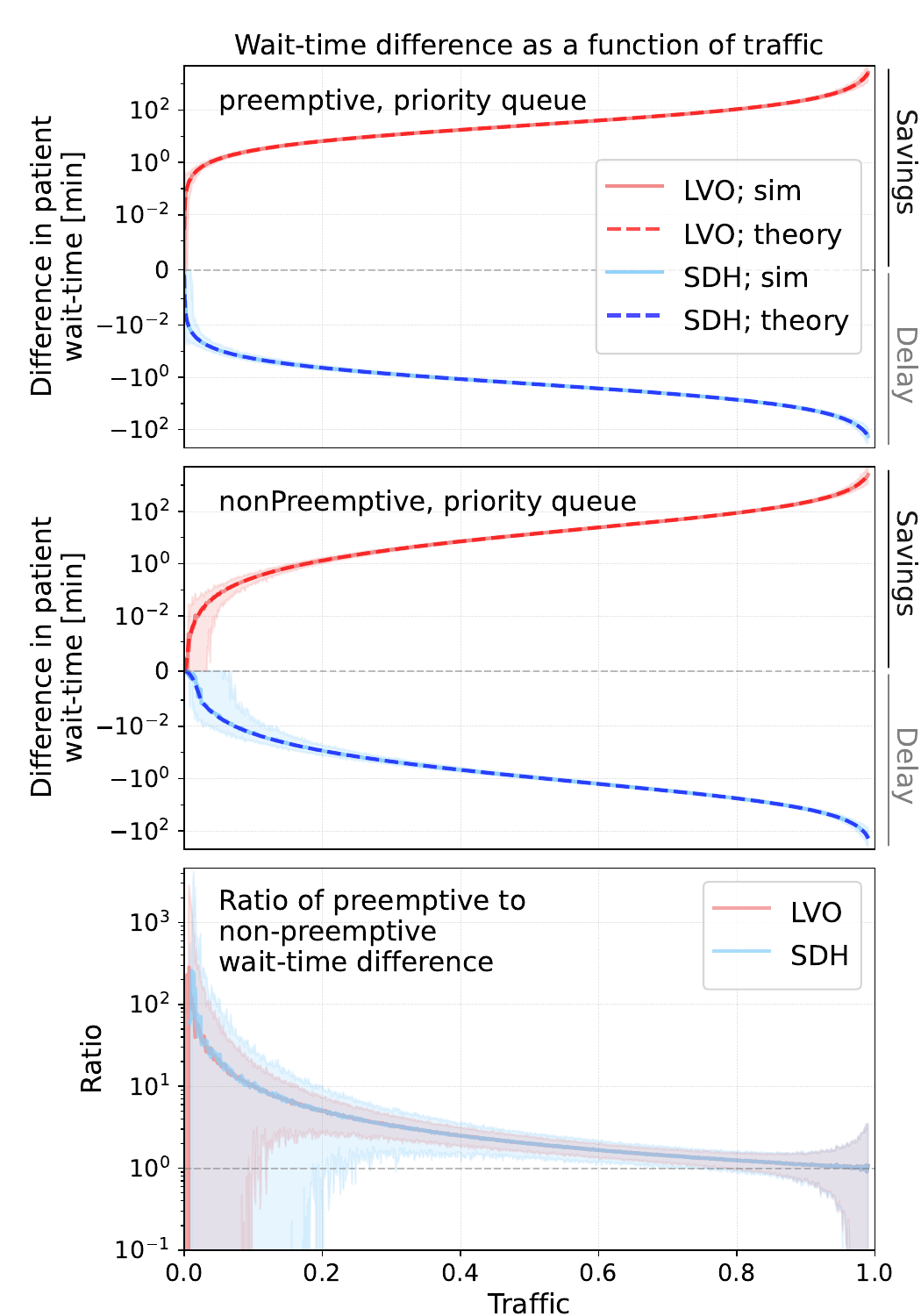}
    &
    \includegraphics[width=0.47\textwidth]{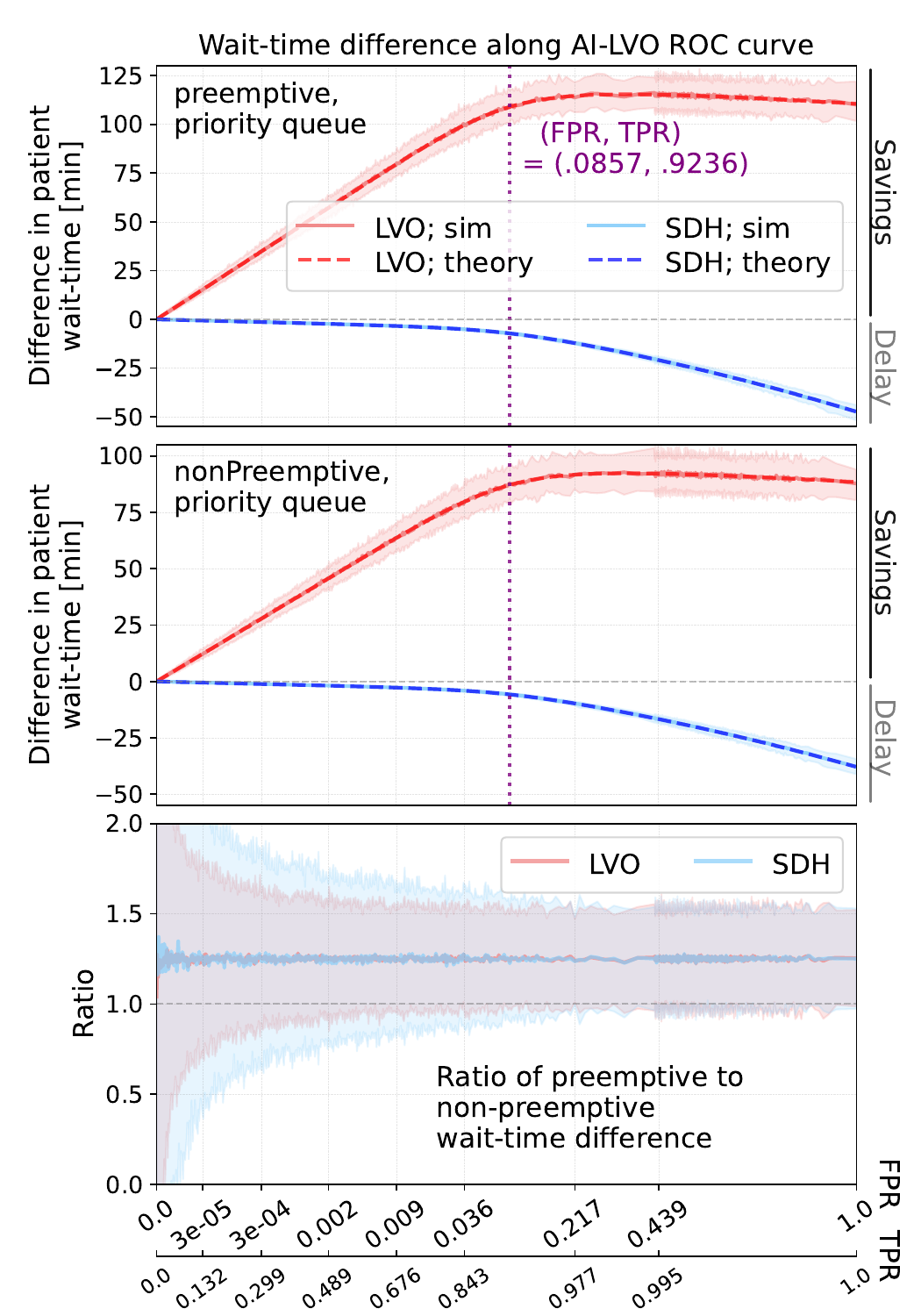}
\end{tabular}    
    \caption{Wait-time results in Experiment 1, assuming a priority protocol, as a function of traffic intensity (left) and AI-LVO performance (right). (Top) Mean wait-time difference for patients with LVO (red) and for those with SDH (blue) in a preemptive setting. Dashed lines represent theoretical results, and solid lines are simulation results.  (Middle) Same for a non-preemptive discipline. (Bottom) Ratio of wait-time difference in a preemptive setting to that in a non-preemptive setting. For the AI-LVO performance plots, vertical purple dotted lines represent the default AI-LVO operating point. The dual x-axis represents pairs of FPR (1-specificity) and TPR (sensitivity) along the AI-LVO ROC curve. SDH = Subdural Hemorrhage; LVO = Large Vessel Occlusion; AI = Artificial Intelligence; FPR = False-Positive Rate; TPR = True-Positive Rate; ROC = Receiver Operating Characteristic. \label{fig:Exp1_priority}}
\end{figure}

With only AI-LVO, the priority workflow is identical to the hierarchical workflow; hence, only results from the priority workflows are shown. Figure~\ref{fig:Exp1_priority} (top 2, left) shows that, for both preemptive and non-preemptive disciplines, the mean wait-time difference diverges as a function of traffic intensity between LVO and SDH patient images. With an AI-LVO to triage LVO patient images, a time-saving on average is expected for this patient subgroup. However, this time-saving benefits is at the expense of the delay risk for SDH patient images. The more LVO time-savings, the more SDH delay. However, the amount of LVO time saved is always greater than that of SDH delay. For both LVO and SDH, large errors are observed at low traffic intensity due to high variability in patient wait-times when simulated patients arrive infrequently. 

In Fig.~\ref{fig:Exp1_priority} (top 2, right), the mean wait-time savings for LVO patients increases sharply as the False Positive Rate (FPR) and True Positive Rate (TPR) along the ROC curve increase. For both preemptive and non-preemptive disciplines, the time-saving reaches a maximum at FPR $\sim$ 0.2; the difference of maximum time-saving between the two disciplines is $\sim$25 minutes. At this point, the average time delay for SDH patients is less than $\sim$15 minutes. With more false-positive LVO patients, LVO time-savings declines from its maximum, while the SDH delay grows drastically to $\sim$50 minutes. 

When comparing results from preemptive and non-preemptive disciplines, we found that wait-time difference is larger in a preemptive discipline across the full traffic intensity spectrum. This is because preemptive discipline saves additional time for incoming LVO cases without waiting for the radiologist to complete the current lower-priority case. The bottom-left plot in Fig.~\ref{fig:Exp1_priority} shows that the ratio in wait-time difference between the two disciplines reaches 1 at high traffic intensity; when the clinic is highly congested, it is most likely that an LVO image is present in the queue, reducing the occurrence of interruption and the distinction between preemptive and non-preemptive disciplines. The large errors at low traffic intensity are primarily because the wait-time in the non-preemptive queue (the denominator of the ratio) is less than 1 minute. In terms of AI-LVO performance, the ratio of preemptive to non-preemptive wait-time remains fairly constant at 1.25 along the ROC curve for both LVO and SDH patient images (Figure \ref{fig:Exp1_priority}, bottom-right plot). These two trends in ratio (reaching to 1 as a function of traffic intensity and staying at a constant with respect to AI performance) are observed across all three experimental workflow scenarios. Therefore, unless otherwise noted, only results from a preemptive discipline are reported in the next two experiment scenarios.

\subsection{Experiment 2: LVO and SDH with two AIs}
\label{results:subsubsec:exp2}
\begin{figure}[h!]
\begin{tabular}{ll}
    \includegraphics[width=0.47\textwidth]{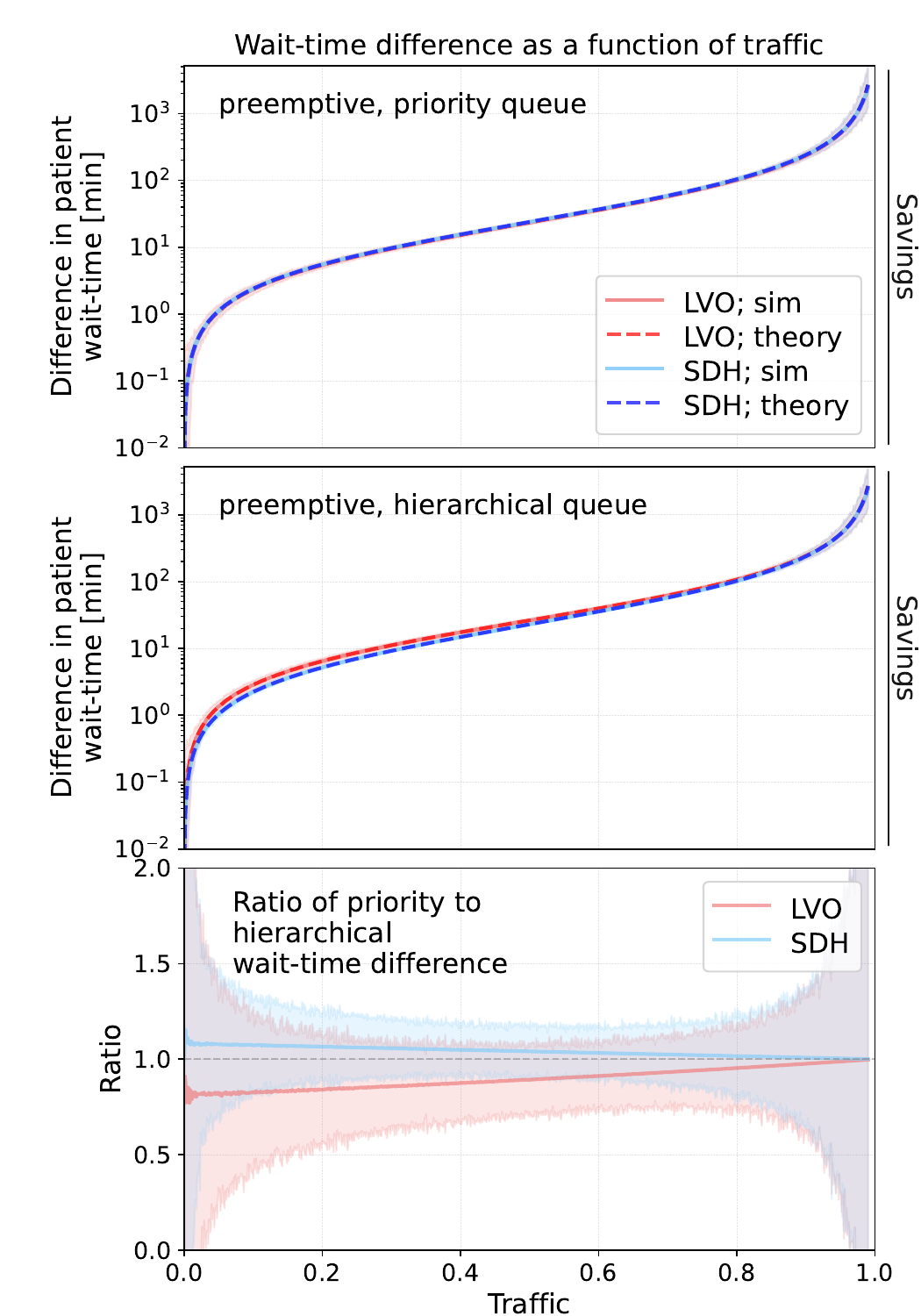}
    &
    \includegraphics[width=0.47\textwidth]{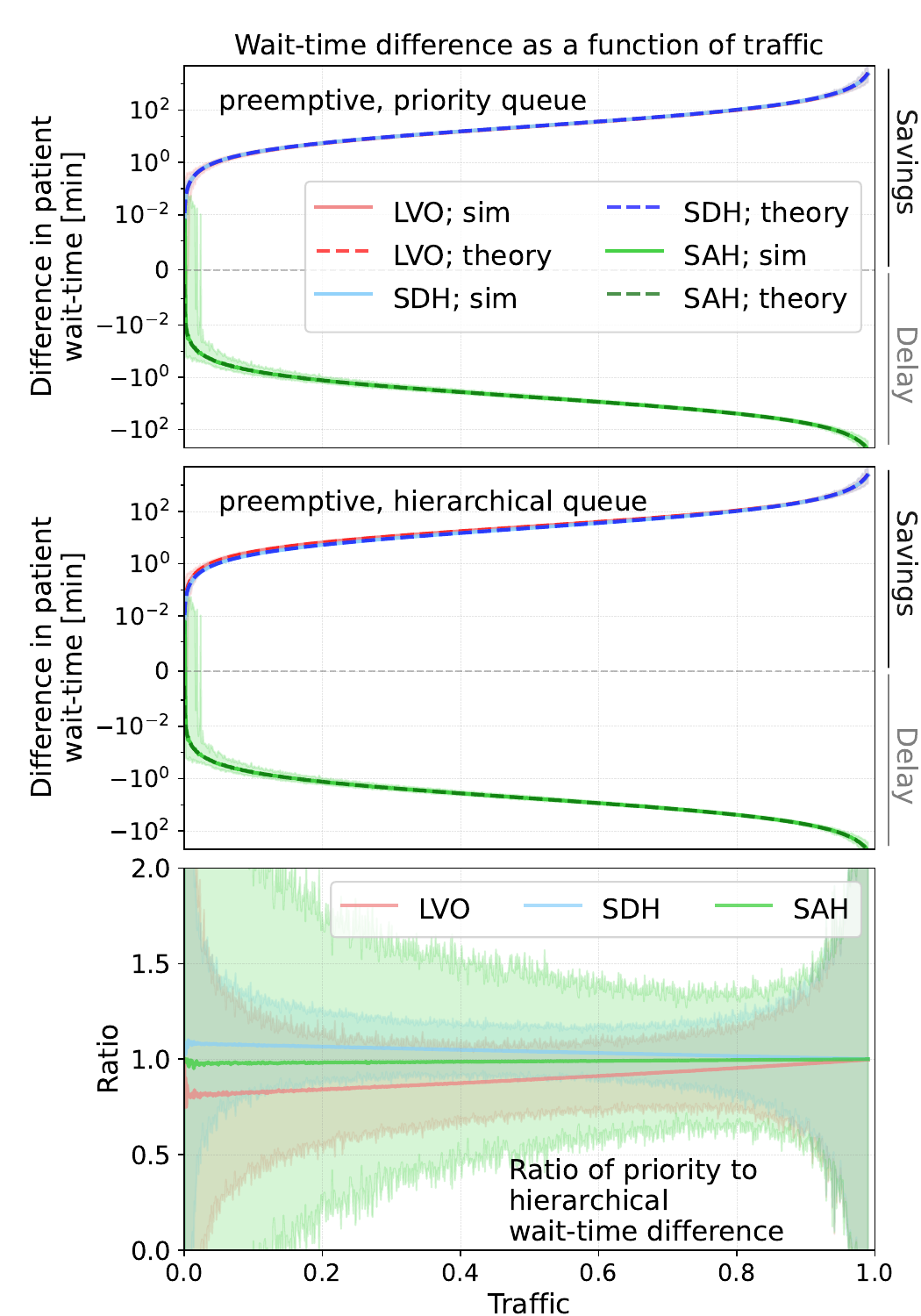}
\end{tabular}
    \caption{Mean wait-time difference in Experiment 2 (left) and 3 (right) as a function of traffic intensity with a preemptive discipline, assuming a priority (top) and hierarchical (middle) protocols. Red lines are time-savings from LVO patient images, and blue lines are that from SDH cases. Experiment 3 also includes the wait-time difference for SAH patient images (green). Dashed lines represent theoretical results, and solid lines are simulation results. (Bottom) Ratio of wait-time difference in the priority workflow to that in the hierarchical workflow. SDH = Subdural Hemorrhage; SAH = Subarachnoid Hemorrhage; LVO = Large Vessel Occlusion; AI = Artificial Intelligence.  \label{fig:Exp2_traffic}}
\end{figure}

As shown in Fig.~\ref{fig:Exp2_traffic} (left), SDH patients are benefited from the addition of AI-SDH into the workflow with limited impact on LVO mean time-savings for both priority and hierarchical protocols. Time-savings for both LVO and SDH patients monotonically increases, suggesting that the use of both AIs (at their default operating points) is beneficial during busy periods in the clinic. The bottom plot suggests that more wait-time is saved for LVO with the hierarchical protocol than the priority protocol because AI-LVO-positive patients are prioritized in the front of AI-SDH-positive patients in a hierarchical setting. As a trade-off, however, the time-savings for SDH patients is slightly diminished compared to the priority setting.

\begin{figure}[h!]
\begin{tabular}{ll}
    \includegraphics[width=0.47\textwidth]{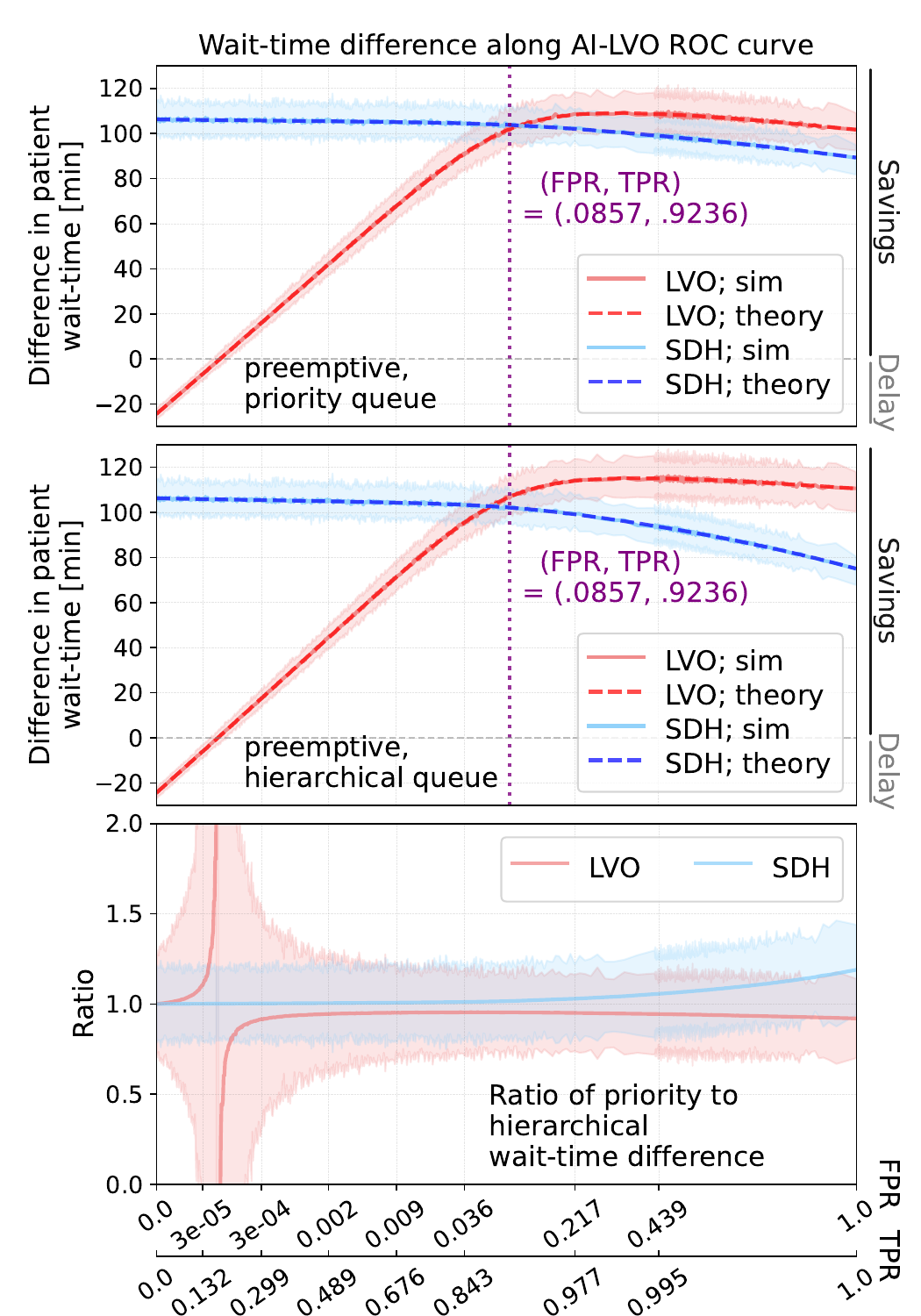}
    &
    \includegraphics[width=0.47\textwidth]{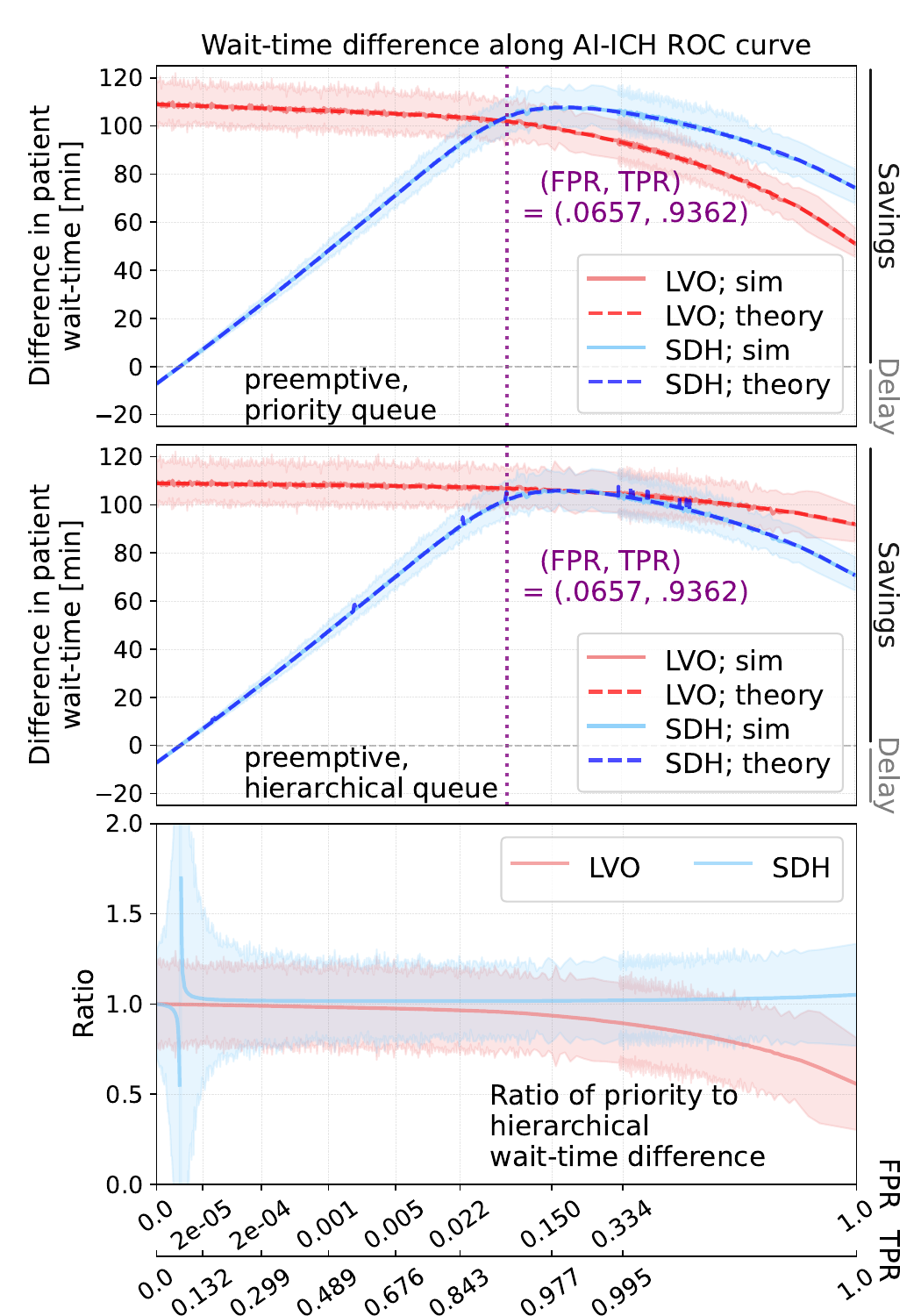}
\end{tabular}    
    \caption{Wait-time results in Experiment 2, assuming a preemptive discipline, as a function of AI-LVO performance when AI-SDH operates at the default operating point (left) and AI-SDH performance when AI-LVO operates at the default Se and Sp (right). (Top) Mean wait-time savings for LVO patients (red) and that for SDH patients (blue) in a priority protocol at the default AI operating points. Dashed lines represent theoretical results, and solid lines are simulation results.  Vertical purple dotted lines represent the default AI operating points for AI-LVO (left) and AI-SDH (right). (Middle) Same for a hierarchical protocol. (Bottom) Ratio of wait-time difference in a priority protocol to that in a hierarchical protocol. Discontinuities occur when the corresponding wait-time difference crosses the zero line. The dual x-axis represents pairs of FPR (1-specificity) and TPR (sensitivity) along the corresponding ROC curve. SDH = Subdural Hemorrhage; LVO = Large Vessel Occlusion; AI = Artificial Intelligence; FPR = False-Positive Rate; TPR = True-Positive Rate; ROC = Receiver Operating Characteristic. \label{fig:Exp2_AI}}
\end{figure}

Figure~\ref{fig:Exp2_AI} shows the dynamic impacts of AI-LVO and AI-SDH operating points on the wait-time difference for LVO and SDH patients. At the default AI-SDH operating point, SDH time-savings slowly diminishes as the AI-LVO-FPR increases along the ROC curve of AI-LVO (top-left plot). This effect is slightly more pronounced in hierarchical protocol (middle-left plot) due to the increasing number of FP patients from AI-LVO that cut in front of the TP SDH patients flagged by AI-SDH. On the other hand, the LVO patients experience a time delay when AI-LVO-FPR = 0 (top/middle-left) since no LVO patients are prioritized, while SDH patients triaged by AI-SDH are reviewed first before these LVO patients. As AI-LVO TPR and FPR increase, LVO patients start to experience more and more time-savings until it reaches a maximum and then gradually decreases. The ratio plot (bottom-left) shows that, at high AI-LVO TPR and FPR, the hierarchical protocol results in slightly more LVO time-savings at the expanse of less ICH time-savings.

The wait-time results in a priority protocol are flipped between SDH and LVO patients when we fix the AI-LVO operating point and vary the AI-ICH performance (Fig.~\ref{fig:Exp2_AI}, top-right). Note that, at the default AI-LVO operating point, the mean time-savings of LVO patients is always above (and, in some cases, the same as) that of SDH in the hierarchical setting regardless of where the AI-SDH operating point is (middle-right). However, this observation does not hold if a lower AI-LVO operating point is chosen, as the middle-left plot shows that SDH time-savings is always above that of LVO in lower AI-LVO TPR and FPR. By comparing the left and right plots in Fig.~\ref{fig:Exp2_AI}, a trade-off between time-savings for patients with the two disease conditions is noticed. By tuning the two AI operating points, an optimal set of operating points can be selected to balance the time-savings for patient images of the two conditions.

\subsection{Experiment 3: LVO, SDH, and SAH with two AIs}

The wait-time difference for the three disease subgroups are shown in Fig.~\ref{fig:Exp2_traffic} (right). The time-savings for LVO and SDH patients are minimally impacted compared to that in Experiment 2 (left of Fig.~\ref{fig:Exp2_traffic}), but the SAH patients are significantly delayed due to the use of AI-LVO and AI-SDH. As shown in the ratio plot (bottom right), the amount of time delay for SAH patients is roughly the same regardless of the protocol because SAH has the lowest priority among the three conditions; the difference in the protocols only concern the re-ordering among patients with the two prioritized conditions only.

\begin{figure}[h!]
\begin{tabular}{ll}
    \includegraphics[width=0.47\textwidth]{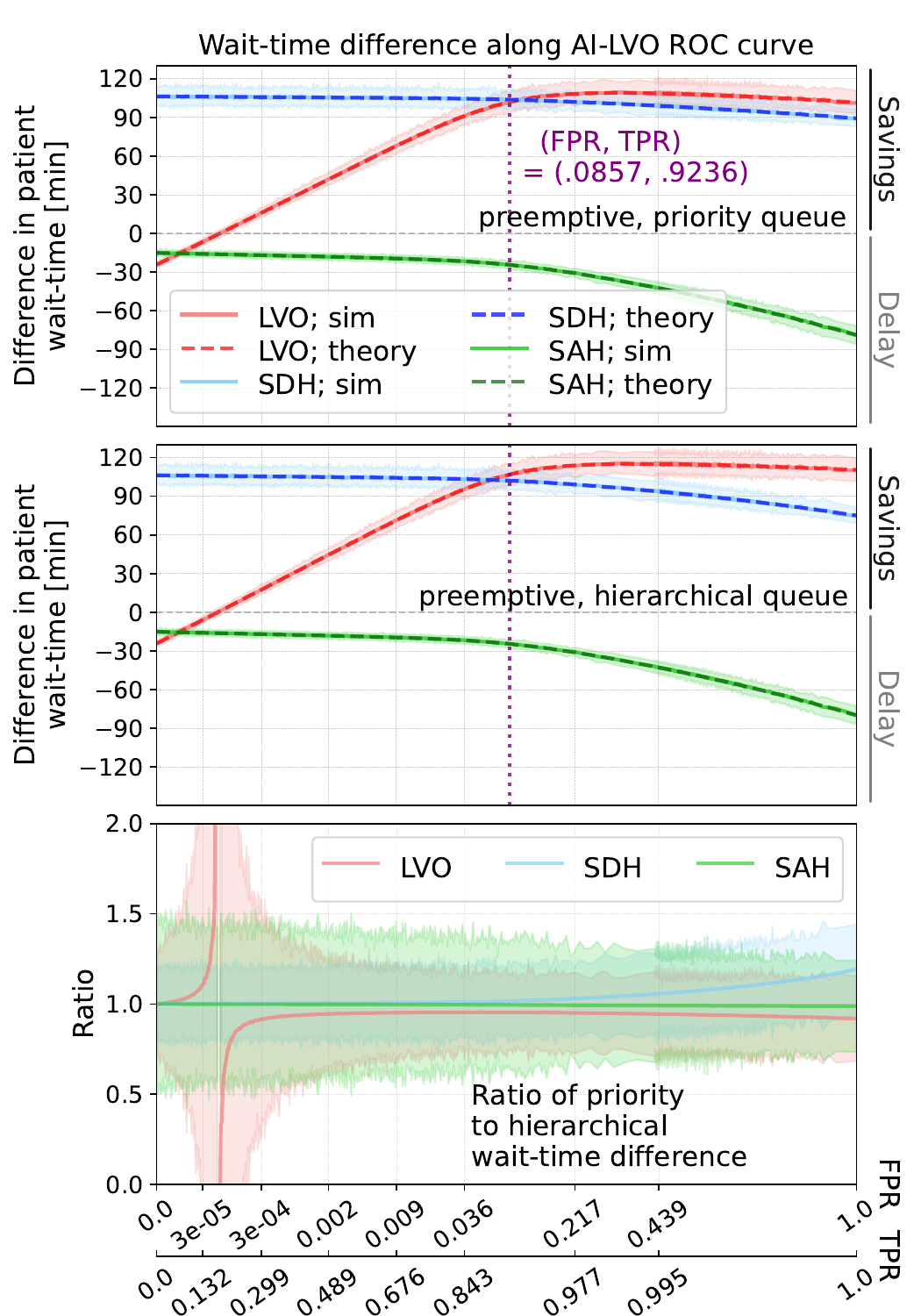}
    &
    \includegraphics[width=0.47\textwidth]{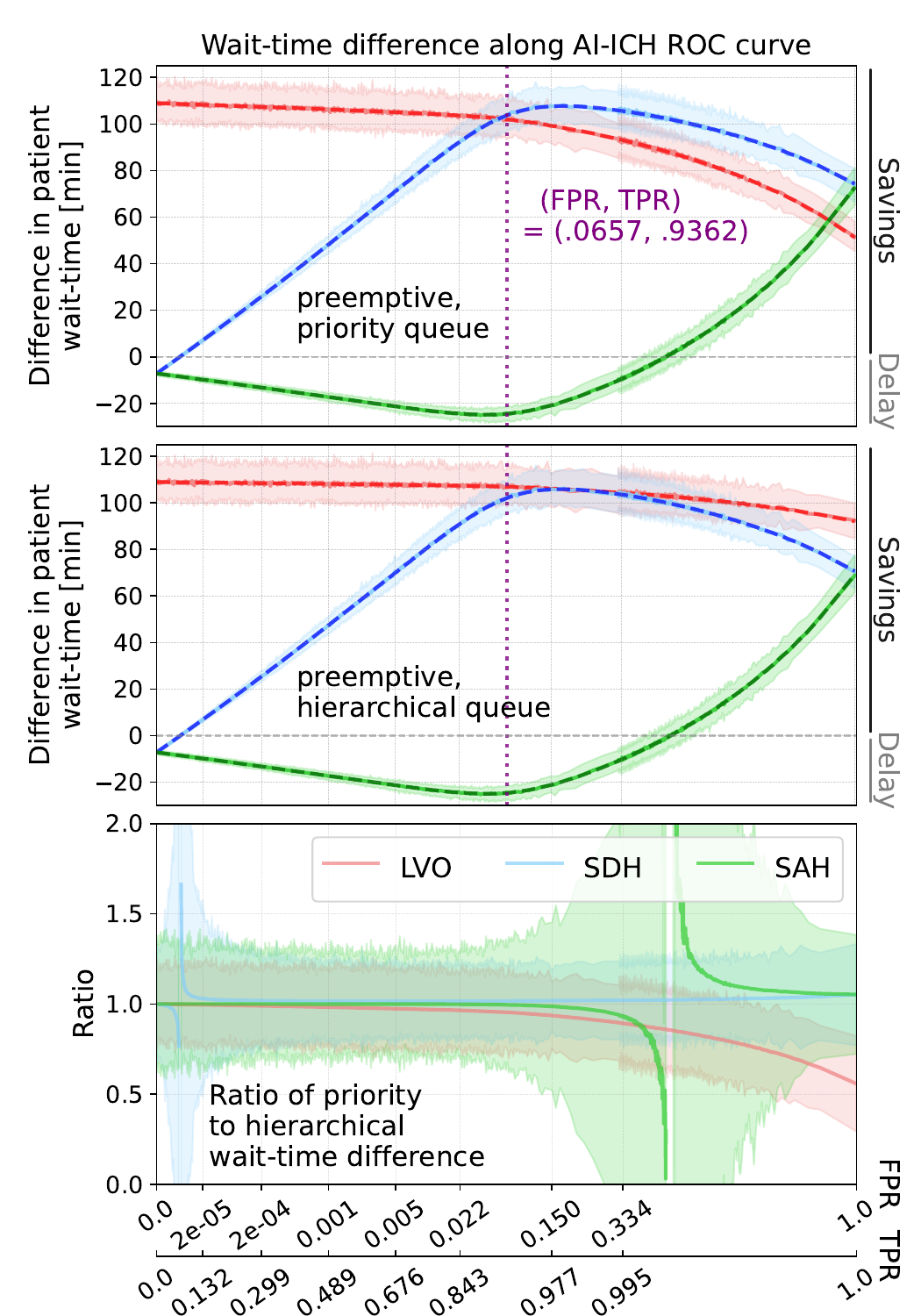}
\end{tabular}    
    \caption{Wait-time results in Experiment 3, assuming a preemptive discipline, as a function of AI-LVO performance (left) and AI-ICH performance (right). (Top) Mean wait-time savings for LVO patients (red), SDH patients (blue), and SAH patients (green) in a priority protocol at the default AI-SDH operating point (left) and at the AI-LVO operating point (right). Dashed lines represent theoretical predictions, and solid lines are simulation results. Vertical purple dotted lines represent the corresponding default AI operating points. (Middle) Same for a hierarchical protocol. (Bottom) Ratio of wait-time difference in a priority protocol to that in a hierarchical protocol. The dual x-axis represents pairs of FPR (1-specificity) and TPR (sensitivity) along the corresponding ROC curve. Discontinuities occur when the corresponding wait-time difference crosses the zero line. SDH = Subdural Hemorrhage; SAH = Subarachnoid Hemorrhage; LVO = Large Vessel Occlusion; AI = Artificial Intelligence; FPR = False-Positive Rate; TPR = True-Positive Rate; ROC = Receiver Operating Characteristic. \label{fig:Exp3_AI}}
\end{figure}

Figure~\ref{fig:Exp3_AI} presents the trade-off between LVO and SDH time-savings and SAH delay as the AI performances vary. The left plots of the figure shows that, at the default AI-SDH operating point, the wait-time difference decreases monotonically for both SDH time-savings and SAH delay; as AI-LVO TPR and FPR increase, prioritizing more CTA images in front, SDH patients experience less time-savings, and SAH patients are more delayed. On the other hand, at a fixed AI-LVO operating point, the SDH and SAH wait-time difference behaves in opposite directions with increasing AI-SDH TPR and FPR. At low TPR and FPR, SDH time-saving increases with increasing SAH time delay when AI-SDH is effective in triaging SDH patients. However, because AI-SDH also processes NCCT images that may have SAH, an increasing AI-SDH FPR mistakenly prioritizes more SAH patients. Hence, SAH time delay turns to time-savings, and the SDH time-savings is diminished.

\subsection{Agreements between Simulation and Theoretical Calculation}
\label{sec:agreement}

Figures~\ref{fig:Exp1_priority}-~\ref{fig:Exp3_AI} show the agreements via REs between simulation and theoretical wait-time difference in Experiment 3 (including all four workflow configurations with preemptive/non-preemptive disciplines and priority/hierarchical protocols) as a function of traffic intensity, disease prevalence, and mean reading time. In addition, we compared simulation and theoretical results in Experiment 4 - the  hypothetical 4-AI, 9-disease workflow.

\begin{figure}[h!]
\centering
\includegraphics[width=1\textwidth]{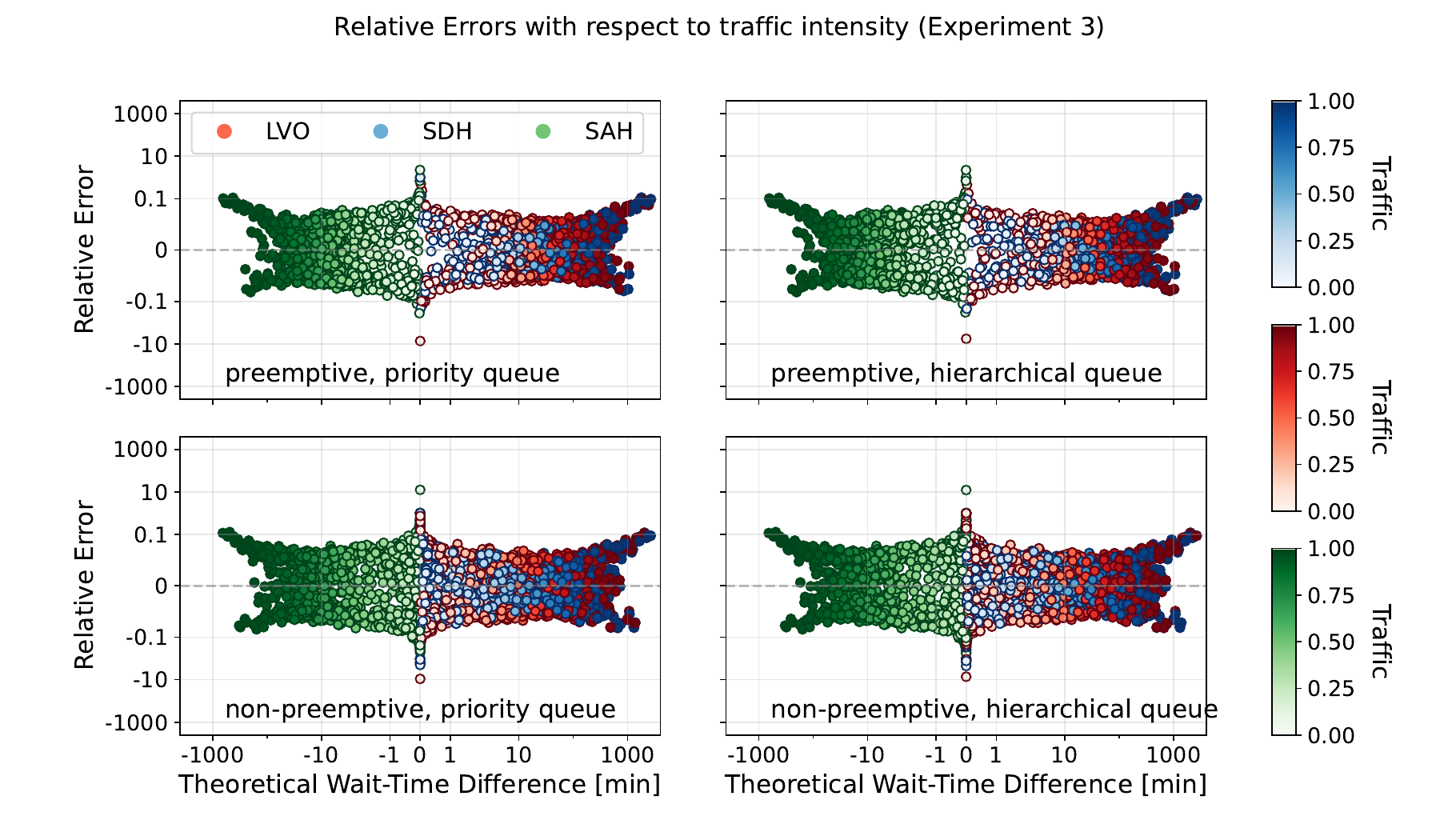}
  \caption{REs between simulation and theoretical calculation in Experiment 3 for preemptive priority (top-left), preemptive hierarchical (top-right), non-preemptive priority (bottom-left), and non-preemptive hierarchical (bottom-right) workflows. REs are presented as a function of theoretical wait-time difference in the x-axis; a positive (negative) value indicates time-savings (delay). Colors represent the disease subgroups: red for LVO patients, blue for SDH patients, and green for SAH patients. Each dot represents a simulation run at a specific traffic intensity indicated by the color gradient; lighter shade means a quieter clinic, while darker shade implies a very congested site. Gray dashed line represents the line of equality i.e. an RE of 0. A vast majority of runs are within an absolute RE of 0.1. SDH = Subdural Hemorrhage; SAH = Subarachnoid Hemorrhage; LVO = Large Vessel Occlusion; RE = Relative Error.\label{fig:agreement_traffic}}
\end{figure}

Figure~\ref{fig:agreement_traffic} shows the REs for traffic intensity ranging from 0 to 1. A vast majority of simulation runs have REs between -0.1 and 0.1. The largest statistical variation happens at very low traffic intensity (i.e. a very quiet clinic without many in-coming patients), where the number of simulated patients with the disease condition of interests is very small. The range of REs is consistent between preemptive (top) and non-preemptive (bottom) disciplines and between priority (left) and hierarchical (right) protocols.

\begin{figure}[h!]
  \centering
  \includegraphics[width=1\textwidth]{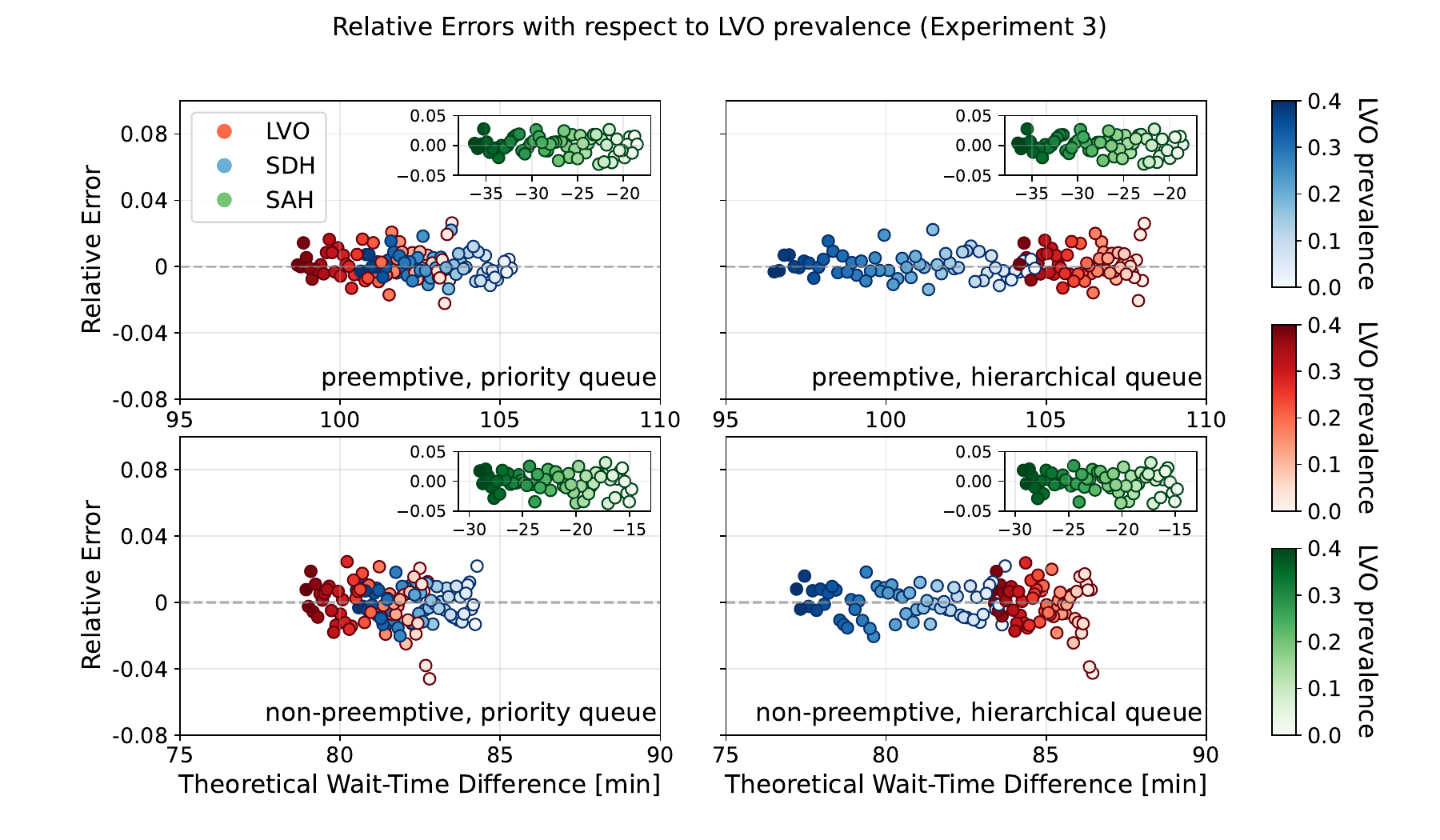}
  \caption{REs between simulation and theoretical calculation in Experiment 3 for preemptive priority (top-left), preemptive hierarchical (top-right), non-preemptive priority (bottom-left), and non-preemptive hierarchical (bottom-right) workflows. REs are presented as a function of theoretical wait-time difference in the x-axis; a positive (negative) value indicates time-savings (delay). Colors represent the disease subgroups: red for LVO patients, blue for SDH patients, and green for SAH patients. SAH results with time delays are presented in separate subplots. Each dot represents a simulation run at a specific LVO prevalence indicated by the color gradient. Gray dashed line represents the line of equality i.e. an RE of 0. A vast majority of runs are within an absolute RE of 0.04. SDH = Subdural Hemorrhage; SAH = Subarachnoid Hemorrhage; LVO = Large Vessel Occlusion; RE = Relative Error.\label{fig:agreement_DiseasePrev}}
\end{figure}

To investigate any disagreement in terms of disease prevalence, simulation runs were generated by varying each disease prevalence from 0 to 0.4, while keeping the other two prevalence rates at their default values. As shown in Fig.~\ref{fig:agreement_DiseasePrev}, REs are within $\pm$0.05 regardless of the disciplines (preemptive vs non-preemptive) and protocols (priority vs hierarchical).  

\begin{figure}[h!]
    \centering
    \includegraphics[width=1\textwidth]{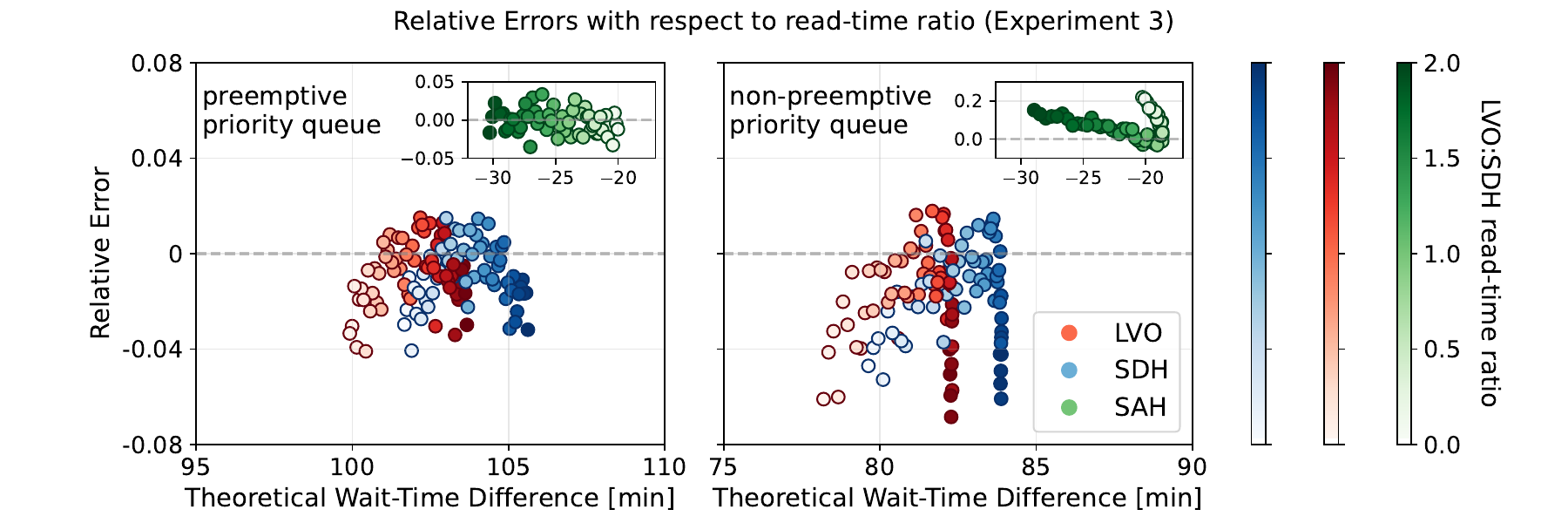}
  \caption{REs between simulation and theoretical calculation in Experiment 3 for preemptive priority (top-left), preemptive hierarchical (top-right), non-preemptive priority (bottom-left), and non-preemptive hierarchical (bottom-right) workflows. REs are presented as a function of theoretical wait-time difference in the x-axis; a positive (negative) value indicates time-savings (delay). Colors represent the disease subgroups: red for LVO patients, blue for SDH patients, and green for SAH patients. SAH results with time delays are presented in separate subplots. Each dot represents a simulation run at a specific ratio between LVO mean read-time and the default 30-min mean read-time for SDH. Indicated by the color gradient, a ratio of 1 refers to a simulation run where the mean read-time are 30 minutes for both LVO and SDH cases. Gray dashed line represents the line of equality i.e. an RE of 0. A vast majority of runs are within an absolute RE of 0.05 with an exception for SAH REs in a non-preemptive priority queue. SDH = Subdural Hemorrhage; SAH = Subarachnoid Hemorrhage; LVO = Large Vessel Occlusion; RE = Relative Error.\label{fig:agreement_readtime}}
\end{figure}
Figure~\ref{fig:agreement_readtime} shows the RE values plotted against wait-time differences when varying the mean read-time for LVO cases from 0 to 60 minutes (i.e. a ratio ranging from 0 to 2 with respect to the default mean read-time of 30 minutes for SDH cases). Only results from the priority protocol are shown because our theoretical approach for the hierarchical protocol is limited to equal mean reading time across all disease conditions. As expected, the RE is the smallest when the read-time ratio is 1 i.e. equal mean read-times for disease condition subgroups. However, as the mean read-times diverge from the default value, an under-estimation in the theoretical wait-time difference is observed for LVO and SDH cases with an RE within $\pm$0.05. SAH cases in non-preemptive priority queue has the largest RE (up to 0.25). Nonetheless, this relatively large RE is still an order of magnitude smaller than the expected wait-time difference; at a 20-minute time delay, an RE of 0.25 implies a difference of 5 minutes between simulation and theoretical calculation, which is still small compared to 20 minutes.

\begin{figure}[h!]
  \centering
    \includegraphics[width=0.85\textwidth]{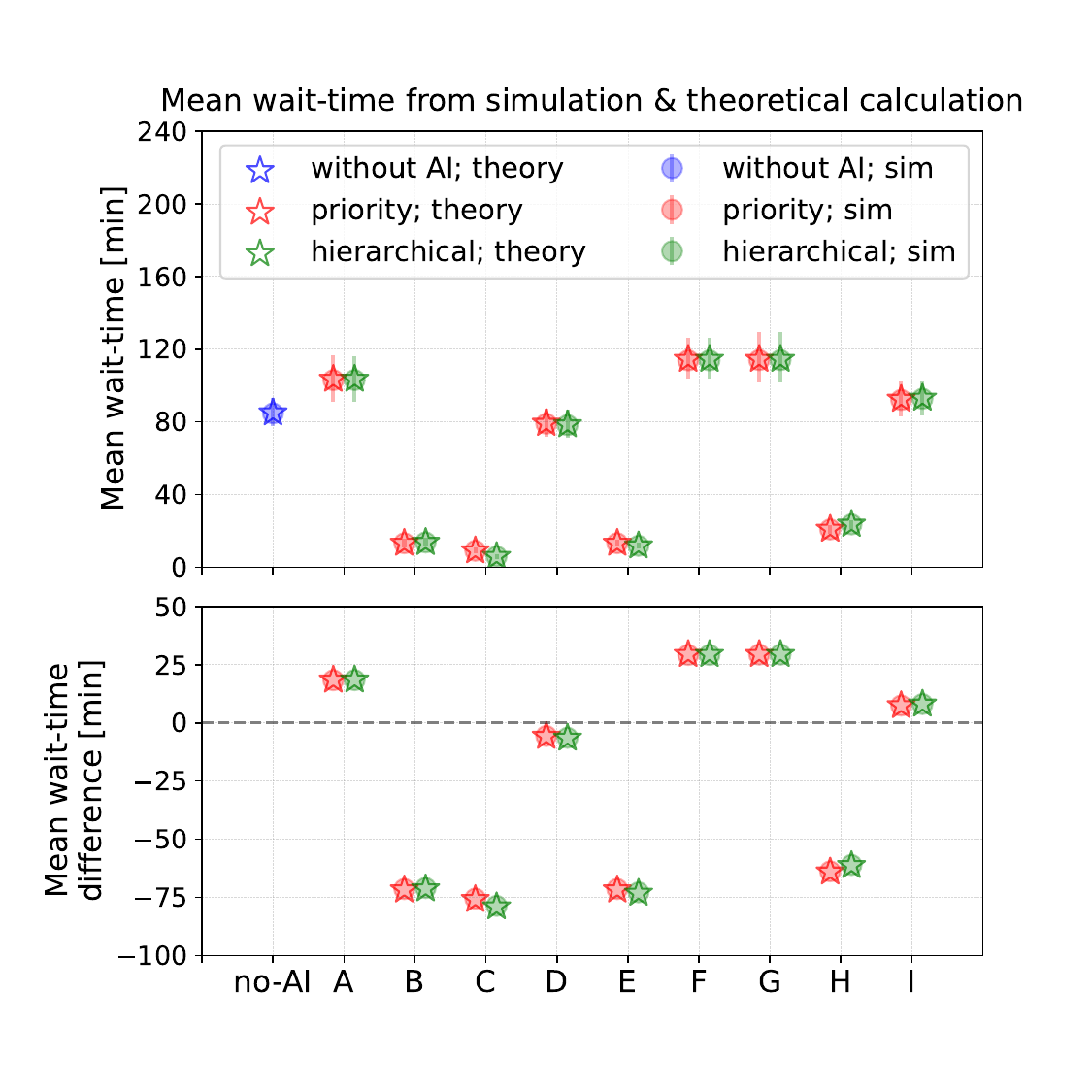}
  \caption{Mean wait-time (top) and wait-time difference (bottom) for all disease condition subgroups under the priority (red) and hierarchical (green) protocols in the 4-AI, 9-disease, preemptive workflow. Blue data points represent mean wait-time from the without-AI, FIFO workflow. Stars represent theoretical predictions, and crosses are simulation results with their 95\% C.I.s. For the mean wait-time different plot, a negative (positive) value indicates time-savings (delay). The dashed gray horizontal line represent zero time difference where no time-savings or delay is observed.  FIFO = First-in-first-out; AI = Artificial Intelligence. \label{fig:exp4}}
\end{figure}

Last, in a complex 4-AI, 9-disease workflow, Fig. \ref{fig:exp4} shows the absolute mean wait-times and mean wait-time differences for the nine disease condition subgroups under both priority and hierarchical protocols in a preemptive setting (see Fig.~\ref{fig:hypoParams} for the workflow parameters). The mean wait-time is significantly shorter for patients with disease conditions at which an AI is deployed to target (i.e. B, C, E, and H). Since two AIs (AI-C and -E), both with relatively high FPRs, are involved in analyzing images in Group 2, patients with condition D (in the same image group) on average experience a small time-saving due to the accidental FP by both AIs. Most importantly, Fig. \ref{fig:exp4} shows that, not only do the wait-time differences agree between simulation and theoretical predictions, the wait-time themselves agree for all subgroups. 

\section{Discussion}

We introduced a \emph{multi-QuCAD}, quantitative framework to estimate the wait-time impacts when multiple AI-triage devices are incorporated into a radiologist workflow. The framework includes both theoretical and simulation approaches to cross check each other. With this framework, several brain CT workflows were simulated using workflow parameters with real-world values from the literature, as well as a complex 4-AI, 9-disease workflow. We found that, under certain values of workflow parameters, incorporating AI-triage devices into a multi-disease workflow can yield wait-time reductions for patients with the conditions where there is an AI to identify that condition. However, it also risks deferring the care of false-negative patients or those deemed lower priority but with time-sensitive conditions where no AI is used to identify them. This trade-off is slightly more pronounced in the preemptive discipline than in the non-preemptive setting and in hierarchical protocol than in priority protocol. Comparison of simulation and theoretical results showed small relative errors compared to the order of magnitude of mean wait-time differences in all workflow scenarios.

AI-triage devices are double-edged swords: they speed up diagnoses of the target conditions at the expense of delaying other conditions that do not have an AI to identify. Since the first AI-triage device was FDA-cleared, questions about its risks of de-prioritizing other time-critical conditions have been raised. While multiple publications have shown clinical evidence in triaging time-sensitive conditions, to the best of our knowledge, no prior work exists to quantify the possible delay in care to de-prioritized patients. Our work has shown that the trade-offs between time-savings and delays depend on the clinical settings (prevalence, reading time, traffic intensity), the AI performances, and the patient sub-populations in the entire reading queue. Analyses focusing on the wait-time impacts to the targeted populations of the AI-triage devices (as in Experiment 2) may be blinded from the time delay risks of patients outside of the targeted population (SAH in Experiment 3). Therefore, patient images outside of the targeted populations of the AI-triage devices should be considered when assessing the full benefit-risk profile of these AI-triage devices.

Real-world clinical workflows exhibit significant complexity and variation across radiology departments, making simulation an essential tool for orchestrating multiple AI-triage devices before implementation. To address this heterogeneity, our approach leverages the power of simulation to provide broad applicability across any imaging modalities, disease conditions, and AI-triage systems. Freely available for download, \emph{multi-QuCAD}, a simulation software verified by queueing theory, offers a risk-free environment to test complex interactions and optimize configurations before costly deployment decisions are made. This simulation-driven approach enables \textit{multi-QuCAD} to estimate wait-time impacts and identify potential bottlenecks when AI-triage devices are implemented in any radiology department, particularly in scenarios where multiple AI-triage devices operate concurrently to prioritize different patient subpopulations within the same reading queue. 

This study is limited in three ways. 
\begin{itemize}
    \item \textit{Lack of clinical confirmation}: While we cross-verified the theoretical and simulation approaches, we acknowledged that clinical confirmation of a multi-AI, multi-disease workflow is not yet available due to difficulties in obtaining the data needed for this study. A radiologist reading list typically contains patient cases with many disease conditions. Identifying all the conditions and collecting their timestamps when they enter the reading list, when the case was opened, and when the radiologist reports are completed, involve many layers of technical challenges beyond our control. However,  we would also like to re-iterate that the existing QuCAD software, which focuses on 1-disease, 1-AI workflow, was clinically validated with real-world data~\cite{ThompsonJACR}, providing some confidence that our proposed methodology--which is an extension of that prior work-- has practical applicability despite the current data limitations.   
    \item \textit{Theoretical limitations}: The theoretical approach is not generalized to all possible workflow settings. The theoretical methods for preemptive hierarchical, non-preemptive priority, and non-preemptive hierarchical workflows are only valid for one radiologist reviewing all cases. The mathematical framework for the preemptive hierarchical workflow is also limited to equal mean reading time between different diseased and non-diseased subgroups. However, these theoretical approaches are still valuable and act as check-points to verify the simulation software which can simulate scenarios with any number of radiologists and unequal mean reading time.
    \item \textit{Independence assumptions}: Both simulation and theoretical approaches assume independence among patient images and disease conditions. In the real world, a patient may be scanned more than once; it could be either a re-scanning of the same imaging modality or an immediate follow-up imaging with a different modality. However, the wait-time difference between the with- and without-AI arms is minimal since the effect on mean wait-time is expected to be the same for both arms, and taking the difference cancels out such effect. Our methodology also does not consider correlations among disease conditions, although in reality the occurrence of one disease condition can be more likely due to the presence of another condition. This correlation depends on the clinical context and can be implemented in the future. 
\end{itemize} 

In conclusion, this work presents \emph{multi-QuCAD}~\cite{ThompsonMultiQuCAD}, a comprehensive quantitative framework that reveals the complex trade-offs inherent in deploying multiple AI-triage devices within radiologist workflows. While AI-triage devices can significantly reduce wait-times for their targeted conditions, they simultaneously delay care for patients that are not prioritized, including those with time-sensitive pathologies that lack AI support. The \textit{multi-QuCAD} simulation framework addresses this critical gap by providing device developers, regulatory scientists, and healthcare providers with a risk-free environment to evaluate these trade-offs before implementation. The \emph{multi-QuCAD} software is publicly available (along with its user manual), empowering healthcare decision-makers to make informed decisions that maximize the benefits of AI-triage technology while minimizing unintended consequences for vulnerable patient populations. 

\vspace{6pt} 

\begin{appendices}

\section*{Appendix}
\setcounter{subsection}{0}
\renewcommand{\thesubsection}{A.\arabic{subsection}}

\subsection{Calculation of probabilities $p_i^{+}$ and $\pi_{\text{set}}^{+}$}
\label{App:probPi}
For $i \in \{1, \cdots, N\}$, we want to calculate $p_i^{+}$, defined as the probability among all images in the reading queue that an image is $a_i$-positive but AI-negative for all higher-priority AIs in the same group. If multiple AIs call a case positive, we only consider the prioritization by the AI that targets the condition with the highest priority level. This aligns with the clinical practice that an image is always prioritized to its highest rank possible in the queue.  

In the following calculation, TN, FN, TP, and FP represent true and false negatives and true and false positives respectively. ND means non-diseased, and $B_i$ corresponds to the disease conditions that belongs to the same image group as condition $i$.  The probability $p_i^{+}$ is given by 
\begin{align*}
p_i^{+} &= \mathbb{P}(a_i^{+}, a_j^{-} \text{ for all } j < i \text{ in } B_i) \\
& = g_i \cdot \sum_{k \in B_i} \mathbb{P}(b_k, a_i^{+}, a_j^{-} \text{ for all } j < i)\\
& = g_i \cdot \Big( \sum_{k < i} \mathbb{P}(b_k) \cdot \mathbb{P}(\text{FN } a_k, \text{ TN } a_j \text{ for } j < i, \text{ FP } a_i, \text{ FP or TN } a_j \text{ for } j > i) \\
& \hspace{2cm} + \sum_{k > i} \mathbb{P}(b_k) \cdot \mathbb{P}(\text{TN } a_j \text{ for } j < i, \text{ FP } a_i, \text{ FN or TP } a_k, \text{ FP or TN } a_j \text{ for } j > i) \\
& \hspace{2cm} + \mathbb{P}( b_i) \cdot \mathbb{P}(\text{TN } a_j \text{ for } j < i, \text{ TP } a_i, \text{ FP or TN } a_j \text{ for } j > i) \\
& \hspace{2cm} + \mathbb{P}(\text{ND}) \cdot \mathbb{P}(\text{TN } a_j \text{ for } j < i, \text{ FP } a_i)\Big)\\
& = g_i \cdot \Big[ \sum_{k < i} \pi_k \Big((1-\mathrm{Se}_k) \cdot \prod_{\substack{j < i \\ j \neq k \\ j \in B_i}}  \mathrm{Sp}_j \cdot (1-\mathrm{Sp}_i)\Big) + \sum_{k > i}  \pi_k \Big( \prod_{\substack{j < i \\ j \in B_i}} \mathrm{Sp}_j \cdot (1-\mathrm{Sp}_i) \Big) \\
& \hspace{2cm} + \pi_i \prod_{\substack{j < i \\ j \in B_i}} \mathrm{Sp}_j \cdot \mathrm{Se}_i + (1-\sum_{j \in B_i} \pi_j) \prod_{\substack{j < i \\ j \in B_i}} \mathrm{Sp}_j \cdot (1-\mathrm{Sp}_i) \Big].
\label{eq:AIpos_probs}
\end{align*}

In \S \ref{subsec:Phierarchical}, we use the notation $\pi_{\{L,H\}}^{+}$ to denote the probability (with respect to all images) that an image is AI-positive from one of the AIs in the set. This is simply
$$\pi_{\{L,H\}}^{+} = \sum_{k: a_k \in \{L,H\}} p_k^{+}.$$

\subsection{Calculation of conditional probabilities for deriving effective arrival and reading rates for diseased and non-diseased groups in priority queues}
\label{App:condprob}

This section concerns the effective arrival rates and reading rates in both preemptive and non-preemptive priority queues. For all conditional probabilities below, if multiple AIs call a case positive, we only consider the prioritization by the AI that targets the condition with the highest priority level. This aligns with the clinical practice that an image is always prioritized to its highest rank possible in the queue. 

First, we provided the mathematical expressions for the conditional probabilities that a patient image has a specific condition $j$ (or non-diseased) \emph{given} that it is $a_i$-positive triaged by the AI targeted at condition $i$ (or AI-negative) subgroups i.e. $\mathbb{P}(\text{diseased } b_j \vert a_i^{+} \text{ subgroup})$ and $\mathbb{P}(\text{non-diseased }\vert a_i^{+} \text{ subgroup})$. Recall that the highest priority class is when $j$ = 1. 

\begin{itemize}
\item If $j > i$, the image with condition $j$ should not be prioritized by the AI $a_i$ which is targeted to identify and triage higher priority cases with condition $i$. This prioritization is a false-positive by $a_i$ and a true-negative for other AIs in the same image group targeting at conditions with a priority level higher than condition $i$. 
\begin{align*}
\mathbb{P}&(\text{diseased } b_j \vert a_i^{+} \text{ subgroup}) = \frac{\mathbb{P}(\text{diseased } b_j, a_i^{+} \text{ subgroup})}{\mathbb{P}(a_i^{+} \text{ subgroup})}\\
& = \frac{1}{p_i^+} \cdot \mathbb{P}(\text{diseased } b_j, a_i^{+} \text{ subgroup})\\
& = \frac{1}{p_i^+} \cdot g_i  \cdot \pi_j  \cdot (1-\mathrm{Sp}_i) \cdot \prod_{\substack{k < i \\ k \in B_i}} \mathrm{Sp}_k 
\end{align*}
Here, $p_i^+$ is the probability among all images that an image is $a_i$-positive but negative for all higher-priority AIs in the same image group, which is calculated in~\ref{App:probPi}. $g_i$ and $\pi_j$ are the group probability of the group that $a_i$ belongs to and the disease prevalence of condition $j$ with respect to images within the group. $B_i$ is a set of disease conditions in the corresponding image group.

\item If $j < i$ , the image with condition $j$ is prioritized by $a_i$ (targeting at a condition $i$ with a lower priority level than $j$). This only happens when $a_j$ fails to prioritize the image. Therefore, the conditional probability depends on both the false-negative rate of $a_j$ and false-positive rate of $a_i$, as well as the true-negative rates of all AIs in the image group targeting at conditions with a priority level higher than condition $i$.
\begin{align*}
\mathbb{P}&(\text{diseased } b_j \vert a_i^{+} \text{ subgroup}) = \frac{\mathbb{P}(\text{diseased } b_j, a_i^{+} \text{ subgroup})}{\mathbb{P}(a_i^{+} \text{ subgroup})}\\
& = \frac{1}{p_i^+} \cdot \mathbb{P}(\text{diseased } b_j, a_i^{+} \text{ subgroup})\\
& = \frac{1}{p_i^+} \cdot g_i  \cdot \pi_j \cdot (1- \mathrm{Se}_j) \cdot (1-\mathrm{Sp}_i) \cdot \prod_{\substack{k < i \\ k \neq j \\ k \in B_i}} \mathrm{Sp}_k 
\end{align*}
\item If $j = i$, the image with condition $j$ is correctly prioritized by $a_j$, which depends on the true-positive rate of $a_j$ and true-negative rates of all AIs in the image group targeting at conditions with a priority level higher than condition $j$.
\begin{align*}
\mathbb{P}&(\text{diseased } b_i \vert a_i^{+} \text{ subgroup}) = \frac{\mathbb{P}(\text{diseased } b_i, a_i^{+} \text{ subgroup})}{\mathbb{P}(a_i^{+} \text{ subgroup})}\\
& = \frac{1}{p_i^+} \cdot \mathbb{P}(\text{diseased } b_i, a_i^{+} \text{ subgroup})\\
& = \frac{1}{p_i^+} \cdot g_i  \cdot \pi_i \cdot \mathrm{Se}_i \cdot \prod_{\substack{k < i \\ k \in B_i}} \mathrm{Sp}_k 
\end{align*}

\item The conditional probability that an image in the same image group as $a_i$ does not have any conditions given that it is labeled as  $a_i$-positive is given by
\begin{align*}
\mathbb{P}(\text{non-diseased} \vert a_i^{+} \text{ subgroup}) = 1 - \sum_j \mathbb{P}(\text{diseased } b_j \vert a_i^{+} \text{ subgroup}).
\end{align*}
\end{itemize} 

Second, we provided the equations for the conditional probabilities that a patient image has a condition $j$ (or non-diseased) \emph{given} that it is labeled as AI-negative by all AIs in the image group i.e. $\mathbb{P}(\text{diseased } b_j \vert a^{-} \text{ subgroup})$ and $\mathbb{P}(\text{non-diseased, group } m_i\vert a^{-} \text{ subgroup})$. 

\begin{itemize}
\item For an image with condition $j$ to be labeled as AI-negative by all AIs, it has to be a false-negative by $a_j$ and true-negatives by other AIs in the same image group. 
\begin{align*} \mathbb{P}&(\text{diseased } b_j \vert a^{-} \text{ subgroup}) = \frac{\mathbb{P}(\text{diseased } b_j, a^{-} \text{ subgroup})}{\mathbb{P}(a^{-} \text{ subgroup})}\\
&= \frac{1}{p^{-}} \cdot g_j \cdot \pi_j \cdot (1-\mathrm{Se}_j) \cdot \prod_{\substack{k \neq j \\ k \in B_j}} \mathrm{Sp}_k 
\end{align*}
Here, $p^-$ is the probability among all images that an image is AI-negative by all AIs in the image group i.e.
\begin{align*} 
p^- = 1 - \sum_j p_j^+.
\end{align*}.

\item  For non-diseased (ND) images in an image group $m_i$, the conditional probability depends on the true-negative rates of all the AIs within the same image group.
\begin{align*}
\mathbb{P}(\text{ND, group } m_i \vert a^{-} \text{ subgroup}) &= \frac{\mathbb{P}(\text{ND, group } m_i, \text{ and AI negative})}{\mathbb{P}(a^{-} \text{ subgroup})} \\
& = \frac{1}{p^{-}} \cdot g_i \cdot (1-\sum_{j \in B_i} \pi_j) \cdot \prod_{j \in B_i} \mathrm{Sp}_j.
\end{align*}

\end{itemize} 

\subsection{Conversion from AI-positive and AI-negative wait-times to diseased and non-diseased wait-times for non-preemptive queues}
\label{App:waittime}

By Bayes' theorem, for any AI-positive subgroup $a_i^{+}$, 
\begin{align*}
\mathbb{P}(a_i^{+} \text{ subgroup} \mid \text{diseased } b_j) &= \frac{\mathbb{P}(\text{diseased } b_j \mid a_i^{+} \text{ subgroup}) \cdot \mathbb{P}(a_i^{+} \text{ subgroup})}{\mathbb{P}(\text{diseased } b_j)} \\
& = \frac{\mathbb{P}(\text{diseased } b_j \mid a_i^{+} \text{ subgroup}) \cdot p_i^+}{g_j \cdot \pi_j}
\end{align*}
Appendix~\ref{App:condprob} provides the mathematical expressions for the conditional probabilities $\mathbb{P}(\text{diseased } b_j \mid a_i^{+} \text{ subgroup})$. By definition, $\mathbb{P}(a_i^{+} \text{ subgroup})$ is $p_i^+$ in Appendix~\ref{App:probPi}, and $\mathbb{P}(\text{diseased } b_j)$ is the prevalence of condition $j$ with respect to all images in the reading queue. Similarly, for the AI-negative subgroup $a^{-}$,
\begin{align*}
\mathbb{P}(a^{-} \text{ subgroup} \mid \text{diseased } b_j) &= \frac{\mathbb{P}(\text{diseased } b_j \mid a^{-} \text{ subgroup}) \cdot \mathbb{P}(a^{-} \text{ subgroup})}{\mathbb{P}(\text{diseased } b_j)} \\
& = \frac{\mathbb{P}(\text{diseased } b_j \mid a^{-} \text{ subgroup}) \cdot p^-}{g_j \cdot \pi_j}
\end{align*}
where $\mathbb{P}(\text{diseased } b_j \mid a^{-} \text{ subgroup})$ and $p^-$ are provided in Appendix~\ref{App:probPi}.

The mean wait-time $\overline{W_j}$ for the images with condition $j$ can therefore be calculated as:
\begin{equation*}
\overline{W_j} = \sum_{i=0}^N \overline{W_{a_i}^+} \cdot \mathbb{P}(a_i^{+} \text{ subgroup} \mid \text{diseased } b_j) + \overline{W^{-}} \cdot \mathbb{P}(a^{-} \text{ subgroup} \mid \text{diseased } b_j).
\end{equation*}

\end{appendices}

\section*{Author Contributions}
Conceptualization, Y.L.E.T.; Methodology, M.M., R.D. and Y.L.E.T.; Software, M.M., R.D. and Y.L.E.T.; Verification, M.M., R.D. and Y.L.E.T.; Formal Analysis, M.M., R.D. and Y.L.E.T.; Writing – Original Draft Preparation, M.M.; Writing – Review \& Editing, M.M., R.D., F.S. and Y.L.E.T.; Visualization, M.M. and Y.L.E.T.; Supervision, F.W. and Y.L.E.T.; Project Administration, F.W. and Y.L.E.T.; Funding Acquisition, F.W.

\section*{Funding}
The authors acknowledge funding by appointments to the Research Participation Program at the Center for Devices and Radiological Health administered by the Oak Ridge Institute for Science and Education through an interagency agreement between the U.S. Department of Energy and the U.S. Food and Drug Administration (FDA).

\section*{Acknowledgments}
We acknowledge the assistance of Dr. Weijie Chen and Dr. Gary Levine for their valuable discussions and feedback during the development of this work. We are also grateful to the Center for Devices and Radiological Health (CDRH) High Performance Computing (HPC) team led by Dr. Kenny Cha for providing the necessary computational resources and technical support.

\section*{Conflicts of Interest}
The authors declare no conflicts of interest.

\bibliographystyle{plain}

\begin{thebibliography}{9}

\bibitem{Raicu}
S.~Raicu \textit{et al.}, 
``Effects of the Queue Discipline on System Performance,''
\emph{Applied Math}, vol. 3, no. 1, pp. 37-48, 2023.

\bibitem{Savage}
T.~Savage \textit{et al.}, 
Prospective Evaluation of Artificial Intelligence Triage of Intracranial Hemorrhage on Noncontrast Head CT Examinations,
\emph{American Journal of Roentgenology}, vol. 223, no. 5., 2023.
doi: \href{https://doi.org/10.2214/AJR.24.31639}{10.2214/AJR.24.31639}

\bibitem{Soun}
X.~Soun \textit{et al.}, 
Impact of an automated large vessel occlusion detection tool on clinical workflow and patient outcomes,
\emph{Frontiers in Neurology}, vol. 14, 2023.  doi: \href{https://doi.org/10.3389/fneur.2023.1179250}{10.3389/fneur.2023.1179250}

\bibitem{Thompson}
Y.L.E.~Thompson \textit{et al.},
Applying queueing theory to evaluate wait-time-savings of triage algorithms, \emph{Queueing Systems},  
vol. 108, no. 3-4, pp.~579–610, 2024.

\bibitem{ThompsonJACR}
Y.L.E.~Thompson \textit{et al.},
Impact of artificial intelligence triage on radiologist report turnaround time: Real-world time savings and insights from model predictions, \emph{J. Am. Coll. Radiol.},  
2025. doi: \href{https://doi.org/10.1016/j.jacr.2025.07.033}{10.1016/j.jacr.2025.07.033}

\bibitem{ThompsonQuCAD}
Y.L.E.~Thompson \textit{et al.},
\emph{QuCAD software},
[Online]. Available: \url{https://github.com/DIDSR/QuCAD/}

\bibitem{ThompsonMultiQuCAD}
Y.L.E.~Thompson \textit{et al.},
\emph{multi-QuCAD software},
[Online]. Available: \url{https://github.com/DIDSR/multi-QuCAD/}

\bibitem{Czap}
A.L. Czap \textit{et al.},
Overview of Imaging Modalities in Stroke, \textit{Neurology}, vol. 97, 2021. doi: \href{https://doi.org/10.1212/WNL.0000000000012794}{10.1212/WNL.0000000000012794} 

\bibitem{Penckofer}
M. Penckofer \textit{et al.},
Neuro-imaging in intracerebral hemorrhage: updates and knowledge gaps, \textit{\textit{Frontiers in neuroscience}}, vol. 18, 2024. doi: \href{https://doi.org/10.3389/fnins.2024.1408288}{10.3389/fnins.2024.1408288} 

\bibitem{Chung}
G.H. Chung \textit{et al.},
The comprehensive comparison of imaging sign from CT angiography and noncontrast CT for predicting intracranial hemorrhage expansion: A comparative study, \textit{Medicine}, vol. 101, 2022. doi: \href{https://doi.org/10.1097/MD.0000000000031914}{10.1097/MD.0000000000031914}

\bibitem{openfda}
U.S. Food and Drug Administration,
\emph{OpenFDA},
[Online]. Available: \url{https://open.fda.gov/apis/device/510k/explore-the-api-with-an-interactive-chart/}

\end{thebibliography}

\end{document}